\documentclass[11pt,reqno]{amsart}
\usepackage{amsfonts,amsmath,amssymb,amscd,mathrsfs}
\usepackage{amsfonts, amsbsy, amsmath, amsthm, amssymb, latexsym, verbatim, enumerate, color}
\usepackage{amsfonts, amsbsy, amsmath, amsthm, amssymb, latexsym, verbatim, enumerate,color}
\usepackage{graphicx}

\definecolor{re}{rgb}{1,0.2,0.2}           % Standard colours red, green, blue
 \definecolor{gr}{rgb}{0,1,0}
 \definecolor{bl}{rgb}{0,0,0.6}
 \definecolor{bl2}{rgb}{0,1,0}

\newtheorem{prop}{\sc Proposition}[section]	%General Theorem etc
\newtheorem{cor}[prop]{\sc Corollary}

\newtheorem{lem}[prop]{\sc Lemma}
\newtheorem{thm}[prop]{\sc Theorem}
\newtheorem{con}[prop]{\sc Conjecture}

\newcommand{\pf}{{\it Proof:\quad}}

\newcommand{\cha}{{\rm char}}	%Polynomials
\newcommand{\main}{{\rm main}}
\newcommand{\gga}{\gamma}

\newcommand{\gl}{\lambda}

\newcommand{\Eig}{{\rm Eig}(A,\mu_{i})}

\newcommand{\lx}{{\rm lex}}
\newcommand{\Id}{{\rm I\,}}
\newcommand{\x}{{\rm x}}
\newcommand{\vv}{{\rm v}}

\newcommand{\e}{{\rm e}}
\newcommand{\f}{{\rm f}}
\newcommand{\ee}{{\rm e}}
\newcommand{\re}{\mathbb{R}}  	%Standard Sets/Fields
\newcommand{\cp}{\mathbb{C}}

\newcommand{\ra}{\mathbb{Q}}
\newcommand{\ka}{\mathbb{K}}
\newcommand{\za}{\mathbb{Z}}

\newcommand{\Gal}{{\rm Gal}}
\newcommand{\I}{{\rm I}}
\newcommand{\dne}{\hfill $\Box$ \vspace{0.cm}}

\date{Version 3 June  2020, compiled \today }
\begin{document}

\title[Walk Matrices]{Unlocking the Walk Matrix of a  Graph}
\author{ F. Liu and J. Siemons}

\address{Fenjin Liu, School of Science, Chang'an University, Xi'an, P.R. China, 710046\newline\newline\phantom{xxxxx}{\rm This work is supported by the National Natural Science Foundation of China (No. 11401044), the Fundamental Research Funds for the Central Universities (No. 300102128201) and Foundation of China Scholarship Council (No. 201706565015).}}

\email{fenjinliu@yahoo.com\newline}

%\address{Fenjin Liu,  School of Science, Chang'an University, Xi'an, P.R. China, 710046, and School of Mathematics and Statistics, Xi'an Jiaotong University, Xi'an, P.R. China, 710049\newline\phantom{xx}  {\rm This work is supported by the National Natural Science Foundation of China (Nos. 11401044, 11471005, 11501050), Postdoctoral Science Foundation of China (No. 2014M560754), Postdoctoral Science Foundation of Shaanxi}}
%\email{fenjinliu@yahoo.com\newline}

%\author{\small $^{a,b,c}$Fenjin Liu\footnote{The corresponding author. E-mail addresses: fenjinliu@yahoo.com, j.siemons@uea.ac.uk\newline\hspace*{10pt}{\scriptsize\FiveStarOpen}This work is supported by the National Natural Science Foundation of China (Nos. 11401044, 11471005, 11501050), Postdoctoral Science Foundation of China (No. 2014M560754), Postdoctoral Science Foundation of Shaanxi %the Fundamental Research Funds for the Central Universities (No. 310812161006)and Foundation of China Scholarship Council (No. 201706565015).}\quad$^c$Johannes Siemons\\\small $^a$School of Science, Chang'an University, Xi'an, P.R. China, 710046\\\small $^b$School of Mathematics and Statistics, Xi'an Jiaotong University, Xi'an, P.R. China, 710049\small $^c$School of Mathematics, University of East Anglia,  Norwich, Norfolk, NR4 7TJ, UK.

\address{Johannes Siemons, School of Mathematics, University of East Anglia, NR47TJ, United Kingdom}
\email{J.Siemons@uea.ac.uk}

%{\small Version 29 October, printed \today}}
%\classno{20C33, 20C40}
\maketitle

{\sc Abstract:}\,   Let $G$ be a graph with vertex set $V=\{v_{1},\dots,v_{n}\}$ and adjacency matrix $A.$ For a subset  $S$ of $V$ let $\e=(x_{1},\,\dots,\,x_{n})^{\tt T}$ be the characteristic vector of $S,$ that is, $x_{\ell}=1$ if $v_{\ell}\in S$ and  $x_{\ell}=0$ otherwise. Then  the $n\times n$ matrix $$W^{S}:=\big[{\rm e},\,A{\rm e},\,A^{2}{\rm e},\dots,A^{n-1}{\rm e}\big]$$\\[-15pt] is the {\it walk matrix} of $G$ for $S.$
This name relates to the fact that in $W^{S}$ the $k^{\rm th}$ entry in the row corresponding to $v_{\ell}$ is the number of walks of length $k-1$ from $v_{\ell}$ to some vertex in $S.$

Since $A$ is symmetric the characteristic vector of $S$ can be written uniquely as a sum of eigenvectors of $A.$ In particular, we may enumerate the distinct eigenvalues $\mu_{1},\dots, \mu_{s}$
of $A$ so that  \\[-8pt]
\begin{eqnarray}\label{SSA}{\rm SD}(S)\!:\,\,\,\,\,\,\,\e&=&\e_{1}+\e_{2}+\dots+\e_{r}\,
\end{eqnarray}  \\[-8pt]
where  $r\leq s$ and $\e_{i}$ is an eigenvector of $A$ of $\mu_{i}$ for all $1\leq i\leq r.$  We refer to (\ref{SSA}) as the {\it spectral decomposition} of $S,$ or more properly, of its characteristic vector.

The key result of this paper is that the walk matrix $W^{S}$ determines the spectral decomposition of $S$ and {\it vice versa.} This holds for any non-empty  set $S$ of vertices of the graph and explicit algorithms which establish this correspondence are given. In particular, we show that the number of distinct eigenvectors that appear in \,(\ref{SSA})\, is equal to the rank of $W^{S}.$

Several theorems can be derived from this result. We show that $W^{S}$ determines the adjacency matrix of $G$ if $W^{S}$ has rank $\geq n-1.$ This theorem is best possible as there are examples of pairs of graphs with the same walk matrix of rank $n-2$ but with different adjacency matrices.

As an immediate corollary we have that if $G$ is a graph with a set $S\subseteq V$ for which its walk matrix
$W^{S}$ has rank $\geq n-1$ then some other graph $G^{*}$ with vertex set $V^{*}$  is isomorphic to $G$ if and only if  there is a set $S^{*}\subseteq V^{*}$ so that $W^{S}$ is the same as $W^{S^{*}}$ up to a reordering of  the rows of $W^{S}$.

From  the work of O'Rourke \& Touri~\cite{orourke1} and Tao \& Vu~\cite{tao}  it is known that
$W^{S}$ has rank $n$ for almost all graphs when  $S=V.$ Therefore in particular,  for almost all graphs $W^{V}$ determines the adjacency matrix of the graph. We conjecture that this remains true more generally when  $S$ is an arbitrary non-empty vertex set of the graph. \\[-15pt] %\!\footnote{{\sc Keywords:}\, Walk Matrix of Graphs, Spectral Decomposition, Isomorphism Problem for Graphs. ~~~~\\[1pt]
%\phantom{iiiiii}{\sc Mathematics Subject Classification:}\, 05C50, 05C75, 05E10}\


\section{\sc Introduction}

Let $G$ be a graph with vertex set $V=\{v_{1},\dots,v_{n}\}$ and adjacency matrix $A.$ If $S$ is a subset of $V$ we let $\e=(x_{1},\,\dots,\,x_{n})^{\tt T}$ be its characteristic vector, that is $x_{\ell}=1$ if $v_{\ell}\in S$ and  $x_{\ell}=0$ otherwise, for $\ell=1,\dots,n.$ Then the  $n\times n$ matrix \begin{eqnarray}\label{WW1}W^{S}:=\big[{\rm e},\,A{\rm e},\,A^{2}{\rm e},\dots,A^{n-1}{\rm e}\big],\end{eqnarray}
formed by the  $A^{i}{\rm e}$ as columns, is the {\it walk matrix} of $G$ for the set $S.$ The term walk matrix refers the fact that the $k^{\rm th}$ entry in the row indexed by $v_{\ell}$ is the number of walks in $G$ of length $k-1$ from $v_{\ell}$ to some vertex in $S.$   %When $S=V$ then $W^{S}$ is the {\it standard walk matrix}\, of $G.$
Walk matrices appeared first in 1978 in Cvetkovi\'c~\cite{cvet1}  for the case $S=V$ and in 2012 in Godsil~\cite{godsil3} for arbitrary non-empty $S\subseteq V.$ A survey about walk matrices and the related topic of main eigenvalues and main eigenvectors can be found in~\cite{row1}. More recently walk matrices  have been studied in spectral graph theory~\cite{Wang-Asc, Wang-Gsc, Wang-Asa} and in particular in connection with the question whether a graph is identified  up to isomorphism by its spectrum~~\cite{liu1}.

\medskip

Next let $\mu_{1},\dots,\mu_{s}$ be the  distinct eigenvalues of $G.$ Since $A$ is symmetric any vector in $\re^{n}$ can be written  uniquely as a linear combination of eigenvectors of $A.$  In particular, when we consider a set $S\subseteq V$ as above, it is clear that its characteristic vector $\e$ can be expressed in this way.   For this purpose we renumber the distinct eigenvalues of $A$ suitably so that we can  write  \\[-8pt]
\begin{eqnarray}\label{SS1}{\rm SD}(S)\!:\,\,\,\,\,\,\,\e&=&\e_{1}+\e_{2}+\dots+\e_{r}\,\text{\,\,with $\e_{i}\neq 0$ \,\,and} \nonumber\\
A\e_{i}&=&\mu_{i}\e_{i}\,\,\,\text{for all }\,\,1\leq i\leq r
\end{eqnarray}  \\[-8pt]
for some $r\leq s.$ We refer to (\ref{SS1}) as the {\it spectral decomposition} of $S,$ or more properly, of its characteristic vector.

\medskip
A key result in this paper shows that the walk matrix for $S$ determines the spectral decomposition of $S,$ and {\it vice versa.} This holds for any graph $G$ and any non-empty set $S$ of vertices of $G.$
From  \,(\ref{SS1})\, we define the $n\times r$ {\it eigenvector matrix} $$E^{S}=\big[\e_{1},\, \e_{2},\,\dots ,\,\e_{r}
\big]$$ and corresponding to the eigenvalues $\mu_{1},\,\mu_{2},\,\dots,\,\mu_{r},$  we define the $r\times n$ {\it eigenvalue  matrix} \\[5pt]
$$M^{S}=
 \left(\begin{array}{ccccc}
 1 & \mu_1& \mu_{1}^{2}\,\,\,\cdots\,\,\, & \mu_1^{n-1} \\
 1& \mu_2&\mu_{2}^{2} \,\,\,\cdots\,\,\, & \mu_2^{n-1} \\
 \vdots & \vdots &\vdots \,\,\,\cdots\,\,\, & \vdots \\
 1 & \mu_r& \mu_{r}^{2} \,\,\,\cdots\,\,\, & \mu_r^{n-1}
 \end{array}\right).
$$\\[5pt]
Most results in this paper depend on the following result, see Theorems~\ref{n4.3}\, and~\ref{n4.6} in Sections 3 and 4:

\medskip
\begin{thm}\label{n1.4} Let $G$ be a graph and  let $S$ be a set of vertices of $G.$ Let $W:=W^{S},$ $E:=E^{S}$ and $M:=M^{S}$  as above. Then $$W=E\cdot M.$$
Furthermore, if $W$ is given then both $E$ and $M$ are determined uniquely. In particular, the number of eigenvectors in the spectral decomposition of $S$ is equal to the rank of $W.$
\end{thm}

The remainder of the paper investigates various applications of this theorem. The first concerns the adjacency matrix of the graph:

\medskip
\begin{thm}\label{n1.1}  Let $G$ be a graph and let $S$ be a set of vertices of $G,$  with walk matrix $W:=W^{S}.$
 Suppose that $W$ has rank $\geq n-1.$ Then  $W$ determines the adjacency matrix of $G.$
\end{thm}

We have explicit formulae   for $A$ in terms of $W^{S}$ when  $W^{S}$ has rank $\geq n-1,$ see Theorems~\ref{n5.01} \,and ~\ref{n5.05}.  The theorem is best possible in general: There are graphs $G,\,G^{*},$ with vertex sets $S$ and $S^{*},$ where $W^{S}=W^{S^{*}}$ have rank $n-2$ while the corresponding adjacency matrices are not equal to each other.

To reorder the vertices of the graph amounts to permuting the rows of the walk matrix. It is therefore natural to bring a walk matrix into a standard form by a permutation of its rows. Here we use the lexicographical ordering of the rows of $W^{S}.$ The lex-ordered version of $W^{S}$ is denoted by $\lx(W^{S}),$ see Section 7. If $S=V$ then $W^{S}$ is the {\it standard walk matrix} of $G.$

\medskip
\begin{thm}\label{n1.2}  Let $G$ and $G^{*}$ be  graphs with standard walk matrices $W$ and $W^{*},$ respectively. Suppose that $W$ has rank $\geq n-1.$  Then $G$ is isomorphic to $G^{*}$ if and only if $\lx(W)=\lx(W^{*}).$
\end{thm}

In Section~6 we consider walk equivalence: two graphs $G$ and $G^{*}$ are walk equivalent to each other if their standard walk matrices are the same, $W^{V}=W^{V^{*}}.$ In Theorem~\ref{n6.01}\, we show that $G$ is walk equivalent to $G^{*}$ if and only if their adjacency matrices restrict to the same map on the space generated by the columns of $W^{V}.$

In Section~7 we discuss further applications and probabilistic results.  O'Rourke and Touri~\cite{orourke1}, based on  the work of Tao and Vu~\cite{tao}, have shown that the standard walk matrix is invertible in almost all graphs, see Theorems~\ref{DD1}~and ~\ref{DD2}.  We therefore have the following  application for the graph isomorphism problem:

\medskip
\begin{thm}\label{n1.3}  For almost all graphs $G$ the following holds: $G$ is isomorphic to some other graph  $G^{*}$ if and only if $\lx(W)=\lx(W^{*})$ where
$W$ and $W^{*}$ are the standard walk matrices of $G$ and $G^{*}$ respectively.
\end{thm}

All graphs in this paper are finite, undirected and without loops or multiple edges.  In Section~2 we state the basics required for the spectral decompositions. In Section~3 we introduce main eigenvectors and main eigenvalues. Sections~4 to~7 contain the material  discussed above. Examples and counterexample of certain graphs discussed in the paper are given as an appendix in Section 8.

\section{Preliminaries}

Let $G$ be a graph on the vertex set $V=\{v_{1},\,v_{2},\dots,v_{n}\}.$ For $u$ and $v$ in $V$ we write $u\sim v$ if $u$ is adjacent to $v.$ The adjacency matrix $A=(a_{ij})$ of $G$ is given by $a_{i,j}=1$ if $v_{i}\sim v_{j}$ and $a_{i,j}=0$ otherwise.  The {\it characteristic polynomial} of $G$ is the characteristic polynomial of $A,$ that is $$\cha_{G}(x)=\det(x \Id-A)
$$ where $\Id$ denotes the $n\times n$ identity matrix. A root $\lambda$ of $\cha_{G}(x)$ is an {\it eigenvalue} of $G$ and the collection  $\lambda_1,\, \lambda_2,\dots,\lambda_n$ of  all roots is the  {\it spectrum} of $G.$ Since $A$ is symmetric all eigenvalues are real. We denote the {\it distinct eigenvalues} of $G$ by  $$\mu_{1},\,\mu_{2},\dots,\mu_{s}$$ for a certain $s\leq n,$ in some arbitrary order. (Later it will be essential to reorder the eigenvalues in particular circumstances.) The {\it minimum polynomial} of $G$ is $${\min}_{G}(x)=(x-\mu_{1})(x-\mu_{2})\cdots (x-\mu_{s}).$$
It is well-known that ${\min}_{G}(x)=f(x)$ is the monic polynomial  of least degree with $f(A)=0.$

The smallest field $\ka$ with $\ra\subseteq \ka\subset \re$ which contains all eigenvalues of $G$ is the {\it splitting field} of $\cha_{G}(x),$ denoted  $\ka=\ra[\mu_{1},\dots,\mu_{s}].$ %Any field automorphisms $\gga\!:\,  \ka\to\ka$ fixes all rational numbers;
The set of all field automorphisms $\gga\!:\,  \ka\to\ka$ which map eigenvalues of $A$ to eigenvalues of $A$
forms the {\it Galois group} of $\cha_{G}(x)$ or
of $G,$ denoted  $\Gal(G).$ We express the action of the field automorphism $\gamma$ by $a\mapsto a^{\gamma}$ for $a\in \ka$ and extend this notation to vectors, matrices and polynomials in the obvious way.   For instance, $A^{\gamma}=(a_{ij}^{\gamma})=(a_{ij})=A.$ We will use the fact that $a\in \ka$ belongs to $\ra$ if and only if $a^{\gamma}=a$ for all $\gamma\in \Gal(G).$ A real number that is the root of an integer polynomial with leading coefficient equal to $1$ is an {\it algebraic integer.} Such a number is rational if and only if it is an ordinary integer. Two algebraic integers are {\it algebraically conjugate} if they are roots of the same irreducible monic integer polynomial.

The monic polynomial $f(x)$ of least degree satisfying $f(A)=0$ is the {\it minimum polynomial} of $G,$ denoted $\min_{G}(x).$ Since $A$ is symmetric we have $${\min}_{G}(x)=(x-\mu_{1})(x-\mu_{2})\cdots (x-\mu_{s})\,.$$ Since $(\min_{G}(A))^{\gamma}=\min^{\gamma}_{G}(A)=0$ for all $\gamma\in\Gal(G)$ it follows  that $\min_{G}(x)=\min^{\gamma}_{G}(x)$ and so this is an integer polynomial. Let \begin{equation}\label{n210}\min\phantom{i\!\!\!}_{G}(x)=f_{1}(x)\cdots f_{\ell}(x)\end{equation}  be factored into irreducible integer polynomials $f_{i}(x).$ Then two eigenvalues $\mu$ and $\mu^{*}$ of $G$ are algebraically conjugate  if and only if they are roots of the same polynomial $f_{i}(x)$ for some $1\leq i\leq \ell.$ From  Galois theory~\cite{stewart} we have

\medskip
\begin{thm}\label{n2.11} The orbits of $\Gal(G)$ on the spectrum of $G$ are the equivalence classes of distinct algebraically conjugate eigenvalues of $G.$ \end{thm}

For each $i$ with  $1\leq i\leq s$ consider the polynomial $$m_{i}(x):=(x-\mu_{1})\cdots(x-\mu_{i-1})(x-\mu_{i+1})\cdots(x-\mu_{s})=(x-\mu_{i})^{-1}\min\!\!\phantom{i}_{G}(x)$$ and define the $n\times n$ matrix
\begin{equation}\label{EQ4}E_{i}\,:=\, \frac{1}{m_{i}(\mu_{i})}\,\,m_{i}(A)\,.
\end{equation} Clearly, each $E_{i}$ is symmetric and its coefficients belong to $\ka.$ %\Fe 1: Here I have changed our Notation. It is more common to denote the minimum idempotents by the letters $E_{i},$ rather than $A_{i}$ as before. So now I use ${\rm Eig}(A,\mu)$ for the eigenspace of $A$ for eigenvalue $\mu.$ (Instead of $E_{i}).$ I hope you agree that this is better. \redne

The following results can be found in many books~\cite{Gantmacher, greub,godsil2, Cvetkovic-IntroductionGraphSpec} on linear algebra.  They are also easy to verify directly from the definition. Usually these results are stated for matrices over the real numbers. However, since the polynomials $m_{i}(x)$ are defined over $\ka$ all matrices and  computations are over $\ka.$ This is essential for us: In particular, Galois automorphisms  act on the $E_{i}$ and all associated quantities over $\ka^{n}.$

\medskip
\begin{lem}\label{n2.12} Let $A$ be the adjacency matrix of $G$ with distinct eigenvalues $\mu_{1},\dots,\mu_{s}.$ For $1\leq i\leq s$ let  $E_{i}$ be as above and  denote the $n\times n$ identity matrix by ${\rm I}.$  Then \\
(i)\,\,\,\,  $E^{2}_{i}\,\,=\,\,E_{i}$\, and $E_{i}E_{j}\,\,=\,\,0$ for all $i\neq j\in \{1,\,\dots,\,s\},$\\
(ii) \,\, ${\rm I}\,\,=\,\,E_{1}\,\,+\,\,E_{2}\,\,+\dots +\,\,E_{s}.$\\
(iii)\,\,\,$A\,\,=\,\,\mu_{1}E_{1}\,\,+\,\,\mu_{2}E_{2}\,\,+\dots +\,\,\mu_{s}E_{s}.$
\end{lem}

The matrices $E_{i}$ are the {\it minimum idempotents} or {\it orthogonal idempotents}
of the graph. From the properties in  (i) and (ii) it is easy to compute all powers of $A$ and so we obtain from \,(iii)\, the Principal Equation  for $A$
\begin{equation}\label{EQ11}(*)\quad A^{k}=\mu_{1}^{k}E_{1}+\mu_{2}^{k}E_{2}+\cdots +\mu_{s}^{k}E_{s}\end{equation} for all $k\geq 0.$ (In fact, one may regard this as the `principal axis equation' of the quadratic form attached to $A.$)

\medskip
Let $x\in \ka^{n}$ and consider $E_{i}\!\cdot\! x$ for some $1\leq i\leq s.$ Since  $$A(E_{i}\!\cdot\! x)=(\mu_{1}E_{1}+\mu_{2}E_{2}+\dots +\mu_{s}E_{s})E_{i}\!\cdot\! x=\mu_{i}E_{i}^{2}\!\cdot\!  x=\mu_{i}E_{i}\!\cdot\! x,$$ using Lemma~\ref{n2.12}(i), it follows that   $E_{i}\!\cdot\! x$ is an eigenvector of $A$ for eigenvalue $\mu_{i},$ provided  that $E_{i}\!\cdot\! x\neq 0.$
From  Lemma~\ref{n2.12}(ii) we conclude that \\[0pt]  $$x=E_{1}\!\cdot\! x+\dots +E_{s}\!\cdot\! x.$$\\[0pt] This is the {\it spectral decomposition} of $x$ into eigenvectors of $A.$
Since $E_{i}^{2}=E_{i}$ we can view this matrix as an orthogonal  projection map$E_{i}\!:\,\ka^{n}\to \ka^{n}.$ Its image  $$\Eig:=\,\{E_{i}\!\cdot\!x\,|\, x\in \ka^{n}\,\}$$ is the eigenspace of $A$ for eigenvalue $\mu_{i}.$ Clearly $E_{i}$  and $\Eig$ determine one another, with many useful interconnections. For instance, the multiplicity of $\mu_{i}$ in the spectrum of $G$ is equal to $\dim(\Eig)={\rm trace}(E_{i})$ and so on.

Next suppose that $x$ is an eigenvector of $A$ for eigenvalue $\mu$ and let $\gamma\in\Gal(G).$  Then $Ax=\mu x$ implies $Ax^{\gamma}=\mu^{\gamma}x^{\gamma}$ since $A^{\gamma}=A.$ Thus $x^{\gamma}$ is an eigenvector for eigenvalue $\mu^{\gamma}.$
This shows that the action of $\Gal(G)$ on the set $\{\mu_{1},\dots,\mu_{s}\}$ of all eigenvalues  extends to an action on the set  $\{{\rm Eig}(A,\mu_{1}),\dots,\, {\rm Eig}(A,\mu_{s})\}$ of all eigenspaces and hence also to an action on  the set   $\{E_{1},\dots,\, E_{s}\}$ of all idempotents of $G.$ Further actions will be discussed  in Section 3. We collect these facts:

\medskip
\begin{thm}[Spectral Decomposition]\label{n2.13}Let  $\mu_{1},\dots,\,\mu_{s}$ be the distinct eigenvalues of the graph $G.$ Let $\ka$ be its splitting field and let ${\rm Eig}(A,\mu_{1}),\dots {\rm Eig}(A,\mu_{s})$ be the eigenspaces of its adjacency matrix.  Then
$$\ka^{n}\,\,=\,\, {\rm Eig}(A,\mu_{1})\,\oplus \dots \oplus\, {\rm Eig}(A,\mu_{s})$$ and for $x\in \ka$ we have the spectral decomposition $$x=E_{1}\!\cdot\! x+\dots +E_{s}\!\cdot\! x \quad\text{with} \quad E_{i}\!\cdot\!x \in {\rm Eig}(A,\mu_{i}).$$
For $i\neq j$ there exists a field automorphism $\gamma\in \Gal(G)$ with $\Eig^{\gamma}={\rm Eig}(A,\mu_{j})$ if and only if $\mu_{i}$ is algebraically conjugate to $\mu_{j}.$ %In particular, algebraically conjugate eigenvalues have the same multiplicity and isomorphic eigenspaces.
\end{thm}

\section{\sc Spectral Decomposition and Walk Matrix, I}

Let $G$ be a graph of order $n$ on  vertex set $V=\{v_{1},\dots,v_{n}\}$ with adjacency matrix $A.$ Denote the  distinct eigenvalues of $G$ by $\mu_{1},\dots,\mu_{s}$ and let $\ka=\ra[\mu_{1},\dots,\mu_{s}]$ be the splitting field of $\cha_{A}(x).$

Throughout  $S\subseteq V$ is  a non-empty set of vertices and $\x^{S}:=(x_{1},\dots,x_{n})^{\tt T}\in\ka^{n}$ denotes its {\it characteristic vector}, that is,   $x_{\ell}=1$ if $v_{\ell}\in S,$ and  $x_{\ell}=0$ otherwise. For convenience we write $\e=\x^{S}$ when the context is clear. Since $A$ is symmetric we can express  $\e$ as a sum of eigenvectors of $A.$ For this purpose we renumber the distinct eigenvalues of $A$ suitably so that we can write $\e$ as  \\[-5pt]
\begin{eqnarray}\label{NN31}{\rm SD}(S)\!:\,\,\,\,\,\,\,\e&=&\e_{1}+\e_{2}+\dots+\e_{r}\,\text{\,\,with $\e_{i}\neq 0$ \,\,and} \nonumber\\
A\e_{i}&=&\mu_{i}\e_{i}\,\,\,\text{for all }\,\,1\leq i\leq r.
\end{eqnarray}  \\[-5pt]
Thus  $\e_{i}$ is an eigenvector of $A$ for eigenvalue $\mu_{i}. $ From the general results in Section~2 we have $\e_{i}=E_{i}\e.$   We call (\ref{NN31}) the {\it spectral decomposition} of $\e=\x^{S}$ or of $S.$  Since $A^{k}\e_{i}=\mu^{k} \e_{i}$ for all $k\geq 0$ we can extend  \,(\ref{NN31})\, to all powers of $A,$ writing
\\[-5pt]
\begin{eqnarray}\label{NN31A}{\rm SD}(S)\!:\,\,\,\,A^{k}\e&=&\mu_{1}^{k}\e_{1}+\mu_{2}^{k}\e_{2}+\dots+\mu_{r}^{k}\e_{r}\,\text{\,\,\,for all}\,\,\,\,k\geq 0. \end{eqnarray}  \\[-5pt]
We are interested in the combinatorial significance  of this decomposition.

We call the $\e_{i}$ and $\mu_{i}$  the {\it main eigenvectors} and {\it main eigenvalues} associated to $S,$ respectively. Comments  and references concerning the notion of main  eigenvectors and eigenvalues are available below in Remark 3.1. It is important to  emphasize that the actual  eigenvalues that appear in (\ref{NN31})  depend on $S,$ in general, while the numbering itself is only a matter of convenience.

To the eigenvalues in the decomposition  \,(\ref{NN31})\, we associate the {\it main  polynomial} \\[-3pt] $${\rm main}_{G}^{S}(x):=(x-\mu_{1})\cdot (x-\mu_{2})\cdots(x-\mu_{r})$$ \\[-3pt] for $S$ in $G.$ We denote the Galois group of ${\rm char}_{G}(x)$ by $\Gal(G).$

 \bigskip
 \begin{lem}\label{NN311} Let $S$ be a set of vertices of $G$ with   spectral decomposition   $\x^{S}=e=\e_{1}+\e_{2}+\dots+\e_{r}. $Then we have:\\[5pt]
 (i)\,\,\, For each $\gamma\in \Gal(G)$ the map $\gamma\!:\,\e_{i}\mapsto \e_{i}^{\gamma}$ is a permutation of the set  $\{\e_{1},\e_{2},\dots,\e_{r}\}$ and the map $\gamma\!:\,\mu_{i}\mapsto \mu_{i}^{\gamma}$ is a permutation of the set  $\{\mu_{1},\mu_{2},\dots,\mu_{r}\}.$\\[5pt]
 (ii)\,\, The main polynomial ${\rm main}_{G}^{S}(x)$ has integer coefficients and divides the minimum polynomial of $G.$ It is the unique monic polynomial $f(x)$ of least degree such that $f(A)(\e)=0.$ (In ring theoretical terms,  $f(x)=\main_{G}^{S}(x)$ is the $A$-annihilator of $\e.)$
 \end{lem}

 \pf (i) We have $A^{\gamma}=A$ and $A\e_{i}^{\gamma}=\mu_{i}^{\gamma}\e_{i}^{\gamma}$ from which we conclude that  $\e_{i}^{\gamma}$ is an eigenvector of $A$ for the eigenvalue $\mu_{i}^{\gamma}.$ Since $\e^{\gamma}=\e$ we have $\e_{1}^{\gamma}+\e_{2}^{\gamma}+\dots+\e_{r}^{\gamma}=\e_{1}+\e_{2}+\dots+\e_{r}$ and hence $\{\e_{1}^{\gamma},\e_{2}^{\gamma},\dots,\e_{r}^{\gamma}\}=\{\e_{1},\e_{2},\dots,\e_{r}\},$ since the decomposition into the eigenvectors of a linear map is unique. For the same reason $\mu_{i}^{\gamma}\in \{\mu_{1},\mu_{2},\dots,\mu_{r}\}$ for each $1\leq i\leq r.$

(ii) From (i) we have  $\{\mu_{1}^{\gamma},\mu_{2}^{\gamma},\dots,\mu_{r}^{\gamma}\}=\{\mu_{1},\mu_{2},\dots,\mu_{r}\}$ and so $\big({\rm main}_{G}^{S}(x)\big)^{\gamma}={\rm main}_{G}^{S}(x)$ for all $\gamma\in \Gal(G).$ Therefore ${\rm main}_{G}^{S}(x)$ is an integer polynomial, with leading coefficient $1.$ Since $(A-\mu_{i} \I)(\e_{i})=0$ we have $\big({\rm main}_{G}^{S}(A)\big)\e=0.$ Since  $\big({\rm min}_{G}(A)\big)\e=0$ it follows that  ${\rm main}_{G}^{S}(x)$ divides ${\rm min}_{G}(x),$ as $\za[x]$ is a principle ideal domain.  It is easy to check that if  $f(x)$ is a proper divisor of ${\rm main}_{G}^{S}(x)$ then $f(A)(\e)\neq 0.$\dne

The lemma imposes significant restrictions on the  main eigenvalues that can appear in the spectral decomposition of a set of vertices.   For instance, the following  is immediate from the lemma:

\medskip
\begin{cor}\label{NN312} Suppose that ${\rm min}_{G}(x)$ is irreducible. Then ${\rm main}^{S}_{G}(x)={\rm min}_{G}(x)$ for every non-empty set  $S\subseteq V.$
\end{cor}

{\sc Remark 3.1:}\, The notion of main eigenvalues and main eigenvectors is due to Cvetkovi\'c,~\cite{cvet1}. In the original context, see also~\cite{row1}, an eigenvector is said to be a {\it main eigenvector} if its eigenspace  contains a vector that is not perpendicular to $\e.$ This definition equivalent to the one we are using here:

 \medskip
\begin{lem}\label{n3.0}  The eigenvalue $\mu_{i}$ is a main eigenvalue for $\e$ if and only if ${\rm Eig}(A,\mu_{i})$ contains a vector that is not perpendicular to $\e.$
\end{lem}

\pf Clearly, $E_{i}\e$ belongs to ${\rm Eig}(A,\mu_{i}),$ and if $E_{i}\e\neq 0$ then $\e^{\tt T}(E_{i}\e)=\e^{\tt T}(E_{i}^{2}\e)=(E_{i}\e)^{\tt T}(E_{i}\e)\neq 0$ by Lemma~\ref{n2.12}. Conversely, suppose that $E_{i}a$  is not perpendicular to $e$ for some $a\in \ka^{n}.$ Thus $0\neq(E_{i}a)^{\tt T}\e=a^{\tt T}(E_{i}\e)$ and so $E_{i}\e\neq 0$ is main. \dne

In \cite{cvet1} and~\cite{row1} the set $S$ is equal to $V$ and so $\e=x^{V}$ is the all 1-vector. We refer to this situation as the {\it standard case.}  The general situation, when $S$ is an arbitrary non-empty set of vertices, was first considered in Godsil~\cite{godsil3}.  %Another common definition in the literature~\cite{cvet1, row1} states that $\mu$ is a main eigenvalue if its eigenspace  contains a vector that is not perpendicular to $\e.$ It is easy to check this definition is equivalent to the above.
Various results related to Lemma~\ref{NN311} for the standard  case can be found in Cvetkovi\'c,~\cite{cvet1}, Teranishi~\cite{ter}, Rowlinson~\cite{row1} and Godsil~\cite{godsil2}.

\medskip
We come to the main topic in this paper: To any set $S$ of vertices of the graph $G$ we can associate a so-called {\it walk matrix} $W^{S}.$  We shall see that this matrix  is very closely related to the spectral decomposition of $S.$

As before let $v_{1},\dots, v_{n}$ be the vertices of $G$ and let $k\geq 0$ be an integer. Then a {\it walk of length} $k$ in $G$ is  a sequence of  vertices $w=(u_0,u_1,\ldots,u_k)$ with $u_{i}\in V$ such that $u_{i-1}\sim u_i$ for $i=1,2,\ldots,\,k.$ (These vertices are not necessarily distinct.) For short, we say that $w$ is a {\it $k$-walk\,} from $u_{0}$ to $u_{k},$ and that theses vertices are the {\it ends} of $w.$  When $A$ denotes  the adjacency matrix of $G$ it is well-known that the $(i,j)$-entry of $A^{k}$ is the total number of $k$-walks from $v_{i}$ to $v_{j},$  see for instance~\cite{godsil2, biggs}. Hence we have:

\medskip
\begin{prop}\label{n4.1} Let $S$ be a subset of $V$ with characteristic vector denoted by $\e=\x^{S}$ and let $0\leq k\leq n.$  Then for all $1\leq j\leq n$ the $j$-th entry of $A^{k}\e$ is the total number of $k$-walks from $v_{j}$ to some vertex in $S.$ \end{prop}

{\sc Definition:}  Let $S$ be a subset of $V$ with characteristic vector  $\e:=\x^{s}.$ Then the {\it walk matrix} of $G$ for $S$ is the $n\times n$ matrix with columns $\e,$ $A\e,$ $A^2\e,\dots,\,A^{n-1}\e,$ thus $$W^{S}:=\big[\e,A\e,A^2\e\dots,\,A^{n-1}\e\big].$$ In particular,  the entries in the $j^{\rm th}$ row of $W^{S}$ are  the number of walks of length $0,\dots, n-1$ from $v_{j}$ ending at a vertex in $S.$  When $S=V$ we refer to $W^{S}$ as the {\it standard walk matrix} of $G.$

For convenience we extend this notation: \,For $0\leq i\leq j$ let $$W^{S}_{[i,j]}:=\big[A^{i}\e,A^{i+1}\e,...\,,A^{j-1}\e,A^{j}\e\big].$$ In particular, $W^{S}=W^{S}_{[0,n-1]},$ and $W^{S}_{[0,0]}=\e$ is the {\it first} column of $W^{S}.$

\bigskip
Returning to the spectral decomposition $\x^{s}=\e_{1}+\e_{2}+\dots+\e_{r}$ for $S$ we define the {\it main eigenvector matrix} $E^{S}.$ Its  columns are the main eigenvectors of $S,$ that is  $$E^{S}:=\big[\e_{1},\e_{2},\dots,\e_{r}\big]. $$
This matrix has size $n\times r.$ We also need certain matrices associated  to the main eigenvalues $\mu_{1},\,\mu_{2},\dots,\mu_{r}$ for $S.$  For $0\leq i\leq j$ denote by $M^{S}_{[i,j]}$ the $r\times (j-i+1)$ matrix  \begin{equation}\label{6a}M^{S}_{[i,j]}=
 \left(\begin{array}{cccc}\mu_1^i & \mu_1^{i+1} & \,\,\,\cdots\,\,\, & \mu_1^j \\\mu_2^i & \mu_2^{i+1} & \cdots & \mu_2^j \\\vdots & \vdots &  & \vdots \\\mu_r^i & \mu_r^{i+1} &  \cdots & \mu_r^j\end{array}\right).
\end{equation}
This matrix is the $[i,j]$-{\it main eigenvalue matrix} for $S.$ % (This definition makes sense for arbitrary $0\leq i\leq j.)$
 Note that   $M^{S}_{[i,j]},$  $\main_{G}^{S}(x)$ and the main eigenvalues for $S$ all determine one another when  $i<j$ or when $i=j$ is odd.   (For $i=j\neq 0$ the sign of the eigenvalue may be determined from other information, in some cases.) The $r\times r$ matrix $$M^{S}:=M^{S}_{[0,r-1]}$$ is the {\it main eigenvalue matrix} for $S$ in $G.$  Recall, these definitions refer to the specific arrangement of the eigenvectors for $\e=\x^{S}$ in the decomposition \,(\ref{NN31}). %The rank of main eigenvalue matrices is easy to determine.

\bigskip
\begin{lem}\label{n4.21} (i)\,\,\, For all $0\leq j$ the matrix $M^{S}_{[0,j]}$ has rank $\min\{r, j+1\}.$\\(ii)\,\,For $0<i<j$ the matrix $M^{S}_{[i,j]}$ has rank $ \min\{r, j-i+1\}$ if $0$ is not a main eigenvalue. If $0$ is a main eigenvalue then  $M^{S}_{[i,j]}$ has rank $ \min\{r, j-i+1\}-1.$
\\ (iii) $M^{S}$ is invertible and $\det(M^{S})^{2}$ is an integer.
\end{lem}

\pf  For (i) and (ii) we use  that $M^{S}_{[i,j]}$ is of Vandemond type and so it is easy to compute the determinant of a suitable square submatrix, having in mind that the main eigenvalues are  pairwise distinct. If $0<i$ and $0$ is an eigenvalue remove a row of zeros from $M^{S}_{[i,j]},$ we leave the details to the reader. For (iii) note that $M=M^{S}$ is invertible by (i). By Lemma~\ref{NN311}\, any element $\gamma\in \Gal(G)$ permutes the rows  of $M$ and so  $\det(M^{\gamma})=\pm \det(M).$ This gives $\det((M^{2})^{\gamma})=\det(M^{2})$ for all $\gamma\in \Gal(G)$ and hence $\det(M^{2})$ is an integer.\dne

\bigskip
\begin{thm}\label{n4.3}Let $S$ be a non-empty set of vertices of $G.$ %and assume ${\rm H}(S).$
Then for all $0\leq i\leq j$ we have $W^{S}_{[i,j]}=E^{S}\cdot M^{S}_{[i,j]}.$ In particular, $W^{S}=E^{S}\cdot M^{S}.$
\end{thm}

We note some consequences of this theorem.

\bigskip
\begin{cor}\label{n4.333}Let $S$ be a non-empty set of vertices of $G$ with walk matrix $W^{S}.$ Let $x^{S}=\e=\e_{1}+\e_{2}+\dots+\e_{r}$ be the spectral decomposition for $S.$  Then
$\e_{1},\e_{2},\dots,\e_{r}$ is a basis for the $\ka$-vector space spanned by the columns of  $W^{S}_{[0,r-1]}.$ In particular, $r={\rm rank}(W^{S}_{[0,r-1]})={\rm rank}(W^{S}).$
\end{cor}

\bigskip
\begin{cor}\label{n4.334} Suppose that $\cha_{G}(x)$ is irreducible. Then\\[5pt]
(i)\,\, $\cha_{G}(x)={\rm main}^{S}_{G}(x)$ for all $S\subseteq V,$ and
\\[5pt]
(ii)\, $|V|={\rm rank}(W^{S})$  for all $S\subseteq V.$
\end{cor}

{\it Proof of Theorem~\ref{n4.3}:} Fix some $k$ with $i\leq k\leq j.$ Then $A^{k}\e=\mu_{1}^{k}\ee_{1}+\dots+\mu_{r}^{k}\ee_{r}$ by (\ref{NN31A}). This expression is equal to the $(k+1)^{st}$ column on the right hand side of the equation, as required.  \dne

{\it Proof of Corollary~\ref{n4.333}:} Let $X$ be the vector space spanned by  $\e_{1},\e_{2},\dots,\e_{r}$ over $\ka.$ As the $\e_{i}$ are orthogonal to each other they are linearly independent. Let $Y$ be the vector space spanned by  the columns of $W^{S}_{[0,r-1]}$ over $\ka.$ By  Theorem~\ref{n4.3} we have $W^{S}_{[0,r-1]}=E^{S}\cdot M^{S}_{[0,r-1]}$ and so every column of $W^{S}_{[0,r-1]}$ is a linear combination of $\e_{1},\e_{2},\dots,\e_{r}.$
Hence $Y\subseteq X.$ By Lemma~\ref{n4.21} the matrix $M^{S}_{[0,r-1]}$ has rank $r$ and so there is a right-inverse $M^{*}$ Thus $M^{S}_{[0,r-1]}M^{*}=\I_{r}$ and hence $W^{S}_{[0,r-1]}M^{*}=E^{S}.$ This implies that $X\subseteq Y.$   \dne

{\it Proof of Corollary~\ref{n4.334}:} Since $\cha_{G}(x)$ is irreducible we have  $\cha_{G}(x)={\rm main}_{G}^{S}(x)$ for any $S\subseteq V,$ by Lemma~\ref{NN311}. By definition, ${\rm main}_{G}^{S}(x)$ has degree $r$ and so $|V|=r.$ Property (ii) now follows from Corollary~\ref{n4.333}.\dne

{\sc Example 3.1:} The simplest kind of walk matrix occurs in regular graphs.  Here $\e_{1}=(1,1,\dots,1)^{\tt T}$ is an eigenvector for  eigenvalue $\mu_{1}=k$ where $k$ is the valency of the graph. Therefore $\x^{V}=\e=\e_{1}$ is the spectral decomposition for $S=V,$ with main polynomial $x-k$ and  walk matrix $W^{V}=[\e,k\e,k^{2}\e,\dots,k^{n-1}\e].$  Connected regular graphs are characterized by this property.

{\sc Example 3.2:} The graph $G$ in Figure 3.1 on the vertex set $V=\{1,2,3,4\}$ has characteristic polynomial $\cha_{G}(x)={\rm min}_G(x)=(1 + x) (x^3-x^{2}-3x+1)$ where the second  factor is irreducible. The Galois group of $\cha_{G}(x)$ is  ${\rm Sym}(1)\times{\rm Sym}(3),$ fixing the rational root
$-1$ and permuting the irrational roots $-1.48...$, $0.31...$ and $2.17...$ as a symmetric group.
%as symmetric group of degree $3.$
\\[8pt] The standard walk matrix of $G$ is $W=W^{V},$ its main polynomial is ${\rm main}_G^{V}(x)=x^3-x^{2}-3x+1.$

\hspace{1.2cm}\includegraphics[width=3.4cm]{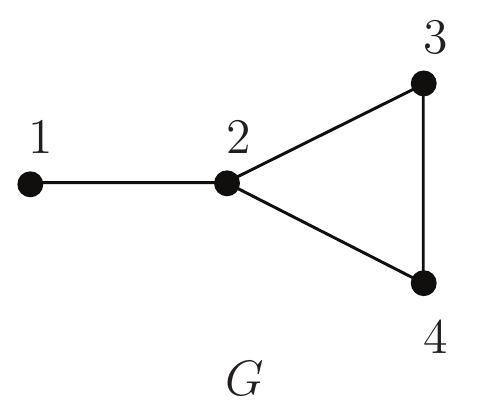}\\[-4.3cm]

\phantom{Z}\hspace{5.3cm}{\small $
  A=\left(
\begin{array}{cccc}
 0 & 1 & 0 & 0 \\
 1 & 0 & 1 & 1 \\
 0 & 1 & 0 & 1 \\
 0 & 1 & 1 & 0 \\
\end{array}\right)$ \hspace{1cm}
$W^{\rm V}=\left(
\begin{array}{cccc}
 1 & 1 & 3 & 5 \\
 1 & 3 & 5 & 13 \\
 1 & 2 & 5 & 10 \\
 1 & 2 & 5 & 10 \\
\end{array}
\right)\,.$}

\bigskip
\centerline{\small \sc Fig 3.1: Graph on 4 Vertices}

\smallskip
We compute the walk matrix and main polynomial for some other subsets.
We have \\[7pt]
{\small $$W^{\{1\}}=\left(
\begin{array}{cccc}
 1 & 0 & 1 & 0 \\
 0 & 1 & 0 & 3 \\
 0 & 0 & 1 & 1 \\
 0 & 0 & 1 & 1 \\
\end{array}\right),\quad
 W^{\{2\}}=\left(
\begin{array}{cccc}
 0 & 1 & 0 & 3 \\
 1 & 0 & 3 & 2 \\
 0 & 1 & 1 & 4 \\
 0 & 1 & 1 & 4 \\
\end{array}
\right),\,  $$\\[2pt]$$W^{\{3\}}=\left(
\begin{array}{cccc}
 0 & 0 & 1 & 1 \\
 0 & 1 & 1 & 4 \\
 1 & 0 & 2 & 2 \\
 0 & 1 & 1 & 3 \\
\end{array}
\right)\, \,\,\text{\,\,and\,\,\,}\,\,\,\,
W^{\{4\}}=\left(
\begin{array}{cccc}
 0 & 0 & 1 & 1 \\
 0 & 1 & 1 & 4 \\
 0 & 1 & 1 & 3 \\
 1 & 0 & 2 & 2 \\
\end{array}
\right).$$}\\[4pt] \!\!\!Their main polynomials are  $\main_{G}^{\{1\}}(x)=\main_{G}^{\{2\}}(x)=(x^3-x^{2}-3x+1),$ both of degree $3={\rm rank}(W^{\{1\}})={\rm rank}(W^{\{2\}}),$ according to Corollary~\ref{n4.333}. For the remaining sets we have   $\main_{G}^{\{3\}}(x)=\main_{G}^{\{ 4\}}(x)=(1+x)(x^3-x^{2}-3x+1),$ of degree  $4={\rm rank}(W^{\{3\}})={\rm rank}(W^{\{4\}}).$ According to Proposition~\ref{n4.2} in the next section, the walk matrix for any set $S\subseteq V$ is of the form $W^{S}=\sum_{i\in S} \,W^{\{i\}}.$ %\Fe 2: This is a simple but pleasing observation: All walk matrices can be determined from the walk matrices of singletons. We could think about the linear span $\mathcal{W}$ of such matrices. Is it an algebra?\redne

{\sc Remark} 3.2: Let $S$ be a set and let $W^{S}=(w_{i,j})$ be its walk matrix. We observe that $W^{S}$ contains information about the subgraph $G[S]$ induced on $S\!:$ Evidently $w_{i,1}=1$ if and only if $i\in S.$ Furthermore, according  to Proposition~\ref{n4.1}, we have that $w_{i,2}$ is the number of neighbours of $v_{i}$ in $S.$ Therefore $(w_{i,2})_{w_{i,1}=1}$ is the degree sequence of  $G[S].$ In particular, the second column of $W^{V}$ is the degree sequence of $G.$

{\sc Remark} 3.3: With regard to the spectral decomposition (\ref{NN31}) of a set $S,$ say $x^{S}=\e_{1}+\dots+\e_{r}$, one may ask if the $\e_{i}$ by themselves already determine the $\mu_{i}.$ However, in general this is not the case, as the following shows. \\[7pt]
Consider graphs $G$ and $G^{*}$ with vertex sets $V$ and $V^{*}$ which have the same splitting field $\ka$ and the same number $s$ of distinct eigenvalues $\mu_{1},\dots ,\mu_{s}$ and $\mu^{*}_{1},\dots ,\mu^{*}_{s},$ respectively. Assume furthermore that the vertex sets can be reordered so that ${\rm Eig}(G,\mu_{i})= {\rm Eig}(G^{*},\mu^{*}_{i})$ for all $1\leq i\leq s.$ In this case we call  $G$ and $G^{*}$ {\it eigenspace equivalent.} \\[7pt]
It is easy to see that every regular graph $G$ is eigenspace equivalent to its complement $\overline{G}.$  In this case $x^{S}=\e_{1},\dots, e_{r}$ is expressed by eigenvectors whose eigenvalues may be those of $G$ or of $\overline{G}.$  We have examples of non-isomorphic eigenspace equivalent graphs which are not of this kind, see Appendix 8.1. It is an open problem to characterize graphs by eigenspace equivalence. %\Fe 3. It would be nice to be able to say that the example in Appendix 8.1 is the smallest non-trivial example, ie $G^{*}$ is not the complement of $G.$ Is it feasible to show this? \dner)

{\sc Remark} 3.4: A similar question occurs for the spectral decomposition $x^{S}=\e_{1}+\dots+ \e_{r}$ when we ask if the $\mu_{i}$ by themselves already determine the $\e_{i},$ up to rearranging vertices. Again, in  general this is not the case, as the following shows.\\[7pt]
Consider two graphs $G$ and $G^{*}$ on the same vertex set $V=V^{*}$ which have the same irreducible characteristic polynomial. (In particular, $G$ and $G^{*}$ are co-spectral.) Let $S=V,$ $\e=\x^{S}$ and consider  $e=\e_{1}+\dots +e_{n}=\e^{*}_{1}+\dots+ \e^{*}_{n},$ see Corollary~\ref{NN312}. Then it is easy to show that $G$ is isomorphic to $G^{*}$
if and only if $\e_{1}=\e^{*}_{1},$ up to a permutation of the entries of the vectors. There are no examples if the order of the graphs is $n<8.$ For $n=8$ an example can be found in Appendix 8.6

\section{\sc Spectral Decomposition and Walk Matrix, II}

Let $S$ be a set of vertices of the graph $G$ with adjacency matrix  $A.$ In the last section we have shown that the spectral decomposition of $S$ determines the walk matrix  for $S.$ In the current section we will show that also the converse is true,  namely that the walk matrix for $S$ determines its spectral decomposition. We start with several general properties of walk matrices.

\bigskip
\begin{prop}\label{n4.2} (i) Let $S\subseteq V$ and  let $0\leq i<j.$  Then $AW^{S}_{[i,j]}=W^{S}_{[i+1,j+1]}.$  \\[3pt]
(ii)  Let $S$ and $T$ be disjoint subsets of $V.$ Then $W^{S}\,+\,W^{T}=W^{S\cup T}.$\\[3pt](iii) Let $S\subseteq V$  and  let $0\leq i<j.$ Then
$$(W^{S}_{[i,j]})^{\tt T}\cdot W^{S}_{[i,j]}=\left(\begin{array}{ccccc}n_{2i} & n_{2i+1} &  \cdots & n_{i+j} \\n_{2i+1} & n_{2i+2} &  \cdots & n_{i+j+1} \\%n_{2i+2} & n_{2i+3} &  \cdots & n_{i+j+2} \\
\vdots & \vdots &  \ddots & \vdots \\n_{i+j} & n_{i+j+1} &  \cdots & n_{2j}\end{array}\right) $$ where $n_{k}$ is the number of all $k$-walks in $G$
with both ends in $S.$ (Note, $(W^{S}_{[i,j]})^{\tt T}$ is the transpose of $W^{S}_{[i,j]}.$)
\end{prop}

\pf The statement (i) follows from the definition of $W^{S}_{[i+1,j+1]}.$ For the remainder let ${\rm s},\,{\rm t},\,{\rm u}$ be the characteristic vectors of $S,\,T$ and $U:=S\cup T.$ We denote the corresponding walk matrices by $W^{S},\,W^{T},\, W^{U}.$ The statement (ii) is immediate since ${\rm s}+{\rm t}={\rm u}.$ To prove (iii) let $S=\{a,\,b,\dots\}$ and let  ${\rm a}$,  ${\rm b},\dots$ be the corresponding characteristic vectors. Then $W^{S}_{[i,j]}= W^{\rm a}_{[i,j]}\,+\,W^{\rm b}_{[i,j]} \,+\, \dots$ by (ii). So we expand $(W^{S}_{[i,j]})^{\tt T}\cdot W^{S}_{[i,j]}$ as a sum of terms of the form $X=(W^{\rm a}_{[i,j]})^{\tt T}\cdot W^{\rm c}_{[i,j]}$ with $c\in S.$

A row of $(W^{\rm a}_{[i,j]})^{\tt T}$ is of the form ${\rm a}^{\tt T}A^{k}$ and a column of   $W^{\rm c}_{[i,j]}$ is of the form $A^{\ell}{\rm c}.$ Therefore the corresponding entry of $X$ is ${\rm a}^{\tt T}A^{k}\cdot A^{\ell}{\rm c}={\rm a}^{\tt T}A^{k+\ell}{\rm c}.$ This is the number of $(k+\ell)$-walks from $a$ to $c.$
\dne

In the following, if $A$ and $B$ are graph quantities, we say that \,`\!$A$ determines $B$'\, if there is an algorithm with input $A$ and output $B$ which is independent of any other property of the graph. For instance, Theorem~\ref{n4.3}\, says that for any vertex set $S$ the matrices $E^{S}$ and $M^{S}$ determine the walk matrix $W^{S}.$ Also the converse is true:

\medskip
\begin{thm}\label{n4.6} Let $S$ be a set of vertices of $G$ with walk matrix $W^{S}$ of rank $r.$  \\[5pt]
(i) \,\, If $r<n$ then $W^{S}_{[0,r]}$ determines $W^{S},$ $M^{S}$ and $E^{S}.$ \\[5pt]
(ii) \, If $r=n$ then $W^{S}_{[0,r-1]}(=W^{S})$ determines $M^{S}$ and $E^{S}.$ \\[8pt]
%In particular, $W^{S}$ determines $M^{S}$ and $E^{S}.$
\end{thm}

\pf  We drop superscripts, writing $W^{S}=W$ etc., where possible.   In  (i) we assume  $r<n.$ Notice that the last column of $W_{[0,r]}$ is a linear combination of the first $r$ columns of $W_{[0,r]}$ by Corollary~\ref{n4.333}. Hence there are rational coefficients $f_{0},\dots,f_{r-1}$ so that $$A^{r}\e=f_{0}\e+f_{1}A\e+\dots+f_{r-1}A^{r-1}\e\,.$$ Hence $f(x)= -f_{0}-f_{1}x-\,\dots\,-f_{r-1}x^{r-1}+x^{r}$ is a rational polynomial with $\big(f(A)\big)(\e)=0.$ It follows from Lemma~\ref{NN311} that $f(x)={\rm main}^{S}_{G}(x).$ Hence we have determined all main eigenvalues for $S$ and hence also $M_{[i,j]}$ for all $0\leq i\leq j.$ We have $W_{[0,r-1]}=E\cdot M_{[0,r-1]}$ by Theorem~\ref{n4.3}. Since $M_{[0,r-1]}$ has an inverse $M^{*}$ by Lemma~\ref{n4.21} it follows that $W_{[0,r-1]}\cdot M^{*}=E.$ Therefore  using Theorem~\ref{n4.3} once more we have determined $W_{[i,j]}$ for all $0\leq i\leq j.$

In (ii) we assume $r=n.$ Here $\cha_{G}(x)={\rm main}_{G}^{S}(x)$ by Lemma~\ref{NN311} and since $A$ satisfies its characteristic equation we have
\begin{eqnarray}\label{2.08} 0&=&\cha_{G}(A)=c_{0}{\rm I}+c_{1}A+\cdots +c_{n-1}A^{n-1}+A^{n}\nonumber\\&=&(A-\mu_{1}\Id)(A-\mu_{2}\Id)\cdots (A-\mu_{n}\Id)\,.\end{eqnarray} Here $c_{n-1}=\mu_{1}+\mu_{2}+\cdots +\mu_{n}={\rm trace}(A)=0$ and so
 \begin{eqnarray}\label{2.08a}-A^{n}&=&c_{0}{\rm I}+c_{1}A+\cdots +c_{n-2}A^{n-2}+c_{n-1}A^{n-1}\nonumber\\&=&c_{0}{\rm I}+c_{1}A+\cdots +c_{n-2}A^{n-2}\,.
\end{eqnarray} It follows that \begin{eqnarray}\label{2.09} -A^{n}\e&=&c_{0}\e+c_{1}A\e+\cdots + c_{n-2}A^{n-2}\e\nonumber\\&=&W_{[0,n-2]}\cdot c^{\tt T}\end{eqnarray}
where $c:=(c_{0},c_{1},...,\,c_{n-2}).$ Since $W=W_{[0,n-1]}$ is given we can compute the walk numbers $$n_{n},\,n_{n+1},\dots,\,n_{2n-2}$$ in Proposition~\ref{n4.2}(iii). (These are $(n-1)$ entries from row $2$ to row $n$ in the last column of  $W_{[0,n-1]}^{\tt T}W_{[0,n-1]}.)$
Therefore, if $w:=(n_{n},\,n_{n+1},\dots,\,n_{2n-2}),$ then we have \begin{equation}\label{2.10}W_{[0,n-2]}^{\tt T}\cdot A^{n}\e=w^{\tt T}\,\end{equation} by Proposition~\ref{n4.2}(iii). Taking \,(\ref{2.09})\, and \,(\ref{2.10})\, together we have \begin{equation}\label{2.11}(W_{[0,n-2]}^{\tt T}\cdot W_{[0,n-2]})\,c^{\tt T}=-w^{\tt T}\,.\end{equation} Since ${\rm rank}(W_{[0,n-1]})=n$ the matrix $W_{[0,n-2]}$ and so  also $$W_{[0,n-2]}^{\tt T}\cdot W_{[0,n-2]}$$ have rank $n-1.$ It follows that $W_{[0,n-2]}^{\tt T}\cdot W_{[0,n-2]}$ is invertible and so $$c^{\tt T}=-(W_{[0,n-2]}^{\tt T}\cdot W_{[0,n-2]})^{-1}w^{\tt T}.$$
This means that we have determined $c$ from $W=W_{[0,n-1]}.$ But these are the coefficients of  the characteristic polynomial and so all eigenvalues are determined. Hence $M$ is determined.  Finally apply Theorem~\ref{n4.3}\, to find $E.$\dne

{\sc Remark} 4.1: If $S$ is a set of vertices then we have seen that the spectral decomposition of $S,$ with the matrices  $E^{S}$ and $M^{S},$ determines its walk matrix $W^{S}$ and vice versa. We have examples of graphs $G,$ $G^{*}$ with $V=S=V^{*}=S^{*}$ and matrices $E\neq E^{*}$ but $M=M^{*}.$ For instance, the two graphs $G$ and $G^{*}$ of order $8$ and labelled No.~79 and No.~80 in Cvetkovi\'c,~\cite{Cvetkovic-A-Table-V6} are co-spectral.  They have the same characteristic and main polynomials for $S=V,$ \begin{equation*}
\begin{split}
&{\rm char}_{G}(x)={\rm char}_{G^{*}}(x)=(x^3-x^2-5 x+1)(x-1)(x+1)^2,\\
&{\rm main}^{V}_{G}(x)={\rm main}_{G^{*}}(x)=x^3-x^2-5x+1.
\end{split}
\end{equation*}
Hence they have the same main eigenvalues. But computation shows that  their main eigenvectors are different. Another pair of this kind are the graphs labelled No.~92 and No.~96 in~\cite{Cvetkovic-A-Table-V6}.  \\[8pt]
{\sc Remark} 4.2: We also have  examples of graphs $G,$ $G^{*}$ with $V=S=V^{*}=S^{*}$ and matrices $E= E^{*}$ but $M\neq M^{*}.$ Trivial examples occur for pairs of regular graphs  of the same order but with different valencies.  For non-regular graphs with the same main eigenvectors but different main eigenvalues see Appendix 8.2.%\Fe 6. A similar question: Are the examples we have for Remark 1\&2 examples of least order? \redne

\section{\sc From Walk Matrix to Adjacency Matrix}

Let $S$ be a non-empty set of vertices of the graph $G.$ Here  we  investigate to what degree $W^{S}$ determines the adjacency matrix of $G.$ The following is a prototype of this problem:

 \bigskip
\begin{thm}\label{n5.01} Let $S$ be a set of vertices of $G$ with walk matrix $W=W^{S}.$ Suppose that $W$ has rank $r=n.$  Then $W$ determines the adjacency matrix of $G.$ \end{thm}

\pf Here $W=W_{[0,n-1]}$ determines $E$ and $M=M_{[0,n-1]}$
by Theorem~\ref{n4.6}. % NOW REFERENCE n4.6(ii) XXXX the matrix
From this we find  $M_{[1,n]}$ and hence $W_{[1,n]}=EM_{[1,n]}$ by Theorem~\ref{n4.3}. By Proposition~\ref{n4.2}(i)\,we have $AW_{[0,n-1]}=W_{[1,n]}.$ As $W$ is invertible we obtain $A= EM_{[1,n]}W^{-1}.$ \dne

%If $r=n-1$ then as before $W=W_{[0,n-1]}$ determines the $n-1$ distinct eigenvalues $\mu_{1},\,\dots,\,\mu_{r}$ with eigenvectors $\e_{1},\,\dots,\,\e_{r}.$ The remaining eigenvalue of $A$ is given by $\lambda=-\mu_{1}-\dots-\mu_{r},$ since ${\rm trace}(A)=0,$ and the remaining eigenvector is the orthogonal complement of $\e_{1},\,\dots,\,\e_{r}.$\dne

One important consequence occurs for graphs with irreducible  characteristic polynomial. Here Corollary~\ref{n4.334} and Theorem~\ref{n5.01} provide the following:

\medskip
\begin{thm}\label{n5.01a} Suppose that the characteristic  polynomial of $G$ is irreducible. Then $W^{S}$ determines the adjacency matrix of $G$ for any non-empty set  $S$ of vertices.  \end{thm}

{\sc Remark} 5.1: Some results about the irreducibility of the characteristic polynomial of a graph are available in~\cite{orourke2 }. For probabilistic results on the rank of the standard walk matrix  see Section~8.

Next we come to the case when the walk matrix has rank $r<n.$ This happens  if and only if there are eigenvectors which are not a scalarly dependant on main eigenvectors.

Let $S$ be  non-empty set of vertices. We denote $\e:=x^{S}$ and $W:=W^{S}. $ We  assume  that $r={\rm rank}(W)<n.$  Let  $\mu_{1},\dots,\,\mu_{r}$ and $\ee_{1},\dots,\,\ee_{r}$ be the main eigenvalues and main eigenvectors that appear in the decomposition \\[-8pt]
\begin{eqnarray}\label{NN31B}{\rm SD}(S)\!:\,\,\,\,\,\,\,\e&=&\e_{1}+\e_{2}+\dots+\e_{r}\,.\end{eqnarray}  \\[-8pt] Let $\gl_{r+1},\dots,\,\gl_{n}$ be the remaining non-main eigenvalues with eigenvectors $$\f_{r+1},\dots,\,\f_{n}\in \ka^{n}.\,$$
To avoid any confusion: it may happen that $\gl_{j}=\mu_{i}$ for some $i,\,j.$ Indeed, this will be the case precisely when the eigenspace for $\mu_{i}$ has dimension $>1.$ In this case we take $\f_{j}$ perpendicular to $\e_{i.}$ It follows that the  $\f_{\ell}$ for $r+1\leq \ell \leq n$ are orthogonal to the columns of $W.$ In addition, we select the  $\f_{\ell}$  to be
 an orthonormal set.
Hence \begin{equation}\label{4.10}\f_{r+1},\dots,\,\f_{n}\,\,\text{is an orthonormal  basis of }\,\,\,{\rm ker}(W^{\tt T}).
\end{equation}\\[-0.3cm]
%The rank  of $W_{[0,r-1]}$ is  $r$ by Theorem~\ref{n4.4}\, and we can determine $W$ according to  Theorem~\ref{n4.6}. By Theorem~\ref{n4.4}\, and we can determine $W$ according to  Theorem~\ref{n4.6}.
Next consider the matrix \begin{equation}\label{4.101}\widehat{W}:=\big[W_{[0,r-1]}\,\big|\, \f_{r+1},\f_{r+2},\dots,\f_{n}\big]\,.\end{equation}  Since $W_{[0,r-1]}$ has rank  $r$ by Corollary~\ref{n4.333}\, we conclude from \,(\ref{4.10})\, that $\widehat{W}$ is invertible; the inverse is given  as follows. Since $W_{[0,r-1]}$ has rank $r$ it follows that
$W_{[0,r-1]}^{\tt T}W_{[0,r-1]}$ is invertible and  so we put $$W^{\dagger}:=(W_{[0,r-1]}^{\tt T} W_{[0,r-1]})^{-1}\cdot W_{[0,r-1]}^{\tt T}\,,$$ a matrix of size $r\times n.$ Next let
\begin{equation}\label{n5.13} \overline{W}:=\left(\begin{array}{c}\,\,W^\dagger \\\f^{\tt T}_{r+1} \\\vdots \\\f^{\tt T}_{n}\end{array}\right)\end{equation} and verify that   $\overline{W}\cdot \widehat{W}={\rm I}_{n}$ by using \,(\ref{4.10}). When we compute $\widehat{W}\cdot \overline{W}={\rm I}_{n}$ we obtain  the equation  \begin{equation}\label{n5.14} \sum_{j=r+1}^{n}\,\f_{j}\cdot \f^{\tt T}_{j}\,=\, {\rm I}_{n}\,-\, W_{[0,r-1]}\cdot W^{\dagger}. \end{equation} The matrix  $\f_{j}\cdot \f_{j}^{\tt T}$ is the projection of $\ka^{n}$ onto the hyperplane with normal $\f_{j}$ and so the sum on the left represents  the eigenspace decomposition of ${\rm ker}(W^{\tt T}).$
From \,(\ref{4.101})\, we have $$A\cdot\widehat{W}
=\big[W_{[1,r]}\,\big|\, \gl_{r+1}\f_{r+1},\gl_{r+2}\f_{r+2},\dots,\gl_{n}\f_{n}\big],$$
see also Proposition~\ref{n4.2}(i), and therefore  \begin{eqnarray}\label{n5.15}A&=&\big[W_{[1,r]}\,\big| \,\gl_{r+1}\f_{r+1},\,\gl_{r+2}\f_{f+2},\dots,\,\gl_{n}\f_{n}\big]\,\cdot \,\widehat{W}^{-1}\nonumber\\[7pt]
A&=&\big[W_{[1,r]}\,\big| \,\gl_{r+1}\f_{r+1},\,\gl_{r+2}\f_{f+2},\dots,\gl_{n}\f_{n}\big]\,\cdot \,\overline{W}\nonumber\\[3pt]
A&=&W_{[1,r]}\cdot W^{\dagger}\,\,+\,\,\sum_{j=r+1}^{n}\gl_{j}(\f_{j}\cdot \f_{j}^{\tt T})\,.
\end{eqnarray} Alternatively, this equation can be derived from \,(\ref{n5.14})\, by multiplying both sides by $A$ and using
 Proposition~\ref{n4.2}(i).\\

\medskip
\begin{thm}\label{n5.04} Let $S$ be a set of vertices of the graph $G$ and let $W=W^{S}$ be its walk matrix. Suppose that $r:={\rm rank}(W)\,\leq n=|V|$ and that $\gl_{r+1},\,\dots,\,\gl_{n}$ are the non-main eigenvalues of $G$ for $S.$  Then the adjacency matrix of $G$ is given by $$A=\,W_{[1,r]}\cdot (W_{[0,r-1]}^{\tt T} W_{[0,r-1]})^{-1}\cdot W_{[0,r-1]}^{\tt T}\,\,\,+\,\,\sum_{j=r+1}^{n}\gl_{j}(\f_{j}\cdot \f_{j}^{\tt T})\,$$ where
  $\f_{r+1},\dots,\f_{n}$ is an orthonormal basis of ${\rm ker}\,W^{\tt T}$ consisting of the non-main eigenvectors of $G,$ that is $A\f_{j}=\gl_{j}\f_{j}$ for $r<j\leq n.$
\end{thm}

\pf For $r=n$ the theorem follows from Theorem~\ref{n5.01} and for $r<n$ the result follows from \,(\ref{n5.15}). \dne

{\sc Definition:}\,  Let $S$ be a set of vertices of the graph $G$ and let $W=W^{S}$ be the corresponding walk matrix for $S.$ Then  $$A_{W}\,:=\,W_{[1,r]}\cdot (W_{[0,r-1]}^{\tt T} W_{[0,r-1]})^{-1}\cdot W_{[0,r-1]}^{\tt T}\,$$ is the {\it $W\!$-restriction} of $A,$ or  {\it $S$-restriction} of $A.$
%\Fe 7. I have now changed from 'core' to 'restriction'. This is better since $A_{W}x=Ax$ for any $x$ in the space generated by the columns of $W.$ So $A_{W}$ in some way is the restriction of $A$ to that space.  Do you have a better suggestion?\redne

We note some properties of this matrix.

\medskip

\begin{prop}\label{n5.041} Let $G$ be a graph on the vertex set $V$ with  adjacency matrix $A.$ Suppose  $S\subseteq V$ and let $A_{W}$ be the $W$-restriction of $A$ for $S.$ Then $A_{W}$ is a symmetric matrix with eigenbasis $\e_{1},\dots,\,\e_{r},\,\f_{r+1},\dots\, \f_{n}$ and eigenvalues $\mu_{1},\dots,\mu_{r},\,0,\dots,\,0.$ In particular,\\[1pt]
(i)\,\,\,\, ${\rm rank}\,A_{W}=r$ if $0$ is not a main eigenvalue for $S$ and ${\rm rank}\,A_{W}=r-1$ otherwise;\\
(ii)\,\, $A\cdot A_{W}=A_{W}\cdot A;$  \\
(iii)\, if $X$ denotes the $\ka$-vector space spanned by the columns of $W$ then  $\{A_{W}\cdot x\,|\,x\in \ka^{n}\}\subseteq X$ with equality if and only if $0$ is not a main eigenvalue for $S.$
\end{prop}

%\Fe 8. We should shorten this proposition, not all of it is used later.\redne

\pf Since $A$ and $f_{j}\cdot f^{\tt T}_{j}$ are symmetric also $A_{W}$ is symmetric. It is straightforward to verify from \,(\ref{n5.15})\, that  $\e_{1},\dots,\,\e_{r},\,\f_{r+1},\dots,\, \f_{n}$ is an eigenbasis for the given eigenvalues. To prove (i) notice that the main eigenvalues are always distinct. Verify (ii) on the given basis. For (iii) notice that  $A_{W}\cdot x$ is a linear combination of the columns of $W_{[1,r]}.$ The remainder follows from (i).\dne

{\sc Remark 5.2:} We give some examples for $W$-restrictions. One extremal situation is  the case $A_{W}=A.$ This happens  if and only if all eigenvalues for $S$ are main. In Section 8 we will see that this holds for almost all graphs when $S=V.$  It  also holds for arbitrary $S$  when the characteristic polynomial of $G$ is irreducible, see~Theorem\ref{n5.01a}.  The other extreme occurs  when  $A_{W}$ has rank $1.$ This happens  when $G$ is regular and $S=V.$ Here, if the degree of $G$ is $k,$ then $A_{[W]}=\frac{k}{n}\,J$ where $J$ is the matrix with all entries equal to $1.$

{\sc Remark 5.3:}  Since the roots of a main polynomial are always simple roots the  number of distinct eigenvalues of the graph gives an upper bound for the rank of $A_{W},$ independently of the set $S$ for which $W=W^{S}.$

We come to  applications of Theorem~\ref{n5.04}  when the rank of a walk matrix of the graph is equal  to $|V|-1$ or $|V|-2.$

\bigskip
\begin{thm}\label{n5.05} Let $G$ be a graph with walk matrix $W$ for some subset $S$ of $V.$ Suppose that ${\rm rank}(W)=|V|-1.$ Then $W$ determines the adjacency matrix of $G.$
\end{thm}

\pf Let $\f_{n}$ be a non-main eigenvector of $G$ with eigenvalues $\gl_{n}.$ Then $\f_{n}\cdot f_{n}^{\tt T}$ is determined by $W$ according to \,(\ref{n5.14}).
By Theorem~\ref{n4.6}\, the main eigenvalues $\mu_{1},\dots,\mu_{n-1}$ are determined by $W$ and so $\gl_{n}=-\mu_{1}-\dots-\mu_{n-1}$ is known.
Now $A$ can be determined by Theorem~\ref{n5.04}. \dne

{\sc Remark 5.4:}  A graph as in this theorem has at least $n-1$ distinct eigenvalues.  In Godsil~\cite{godsil0, godsil3} it is shown that any graph with this property is spectrally vertex reconstructible.

\bigskip
\begin{thm}\label{n5.06} Let $G$ be a graph with walk matrix $W$ for some $S\subseteq V.$ Suppose that ${\rm rank}(W)=|V|-2.$ Assume that (a) $S=V$ or (b) $S\neq \emptyset$ is arbitrary and the number of edges of $G$ is given. Then $W$ determines all  eigenvalues of $G.$ Let $\gl_{n-1}$ and $\gl_{n}$ be the two non-main eigenvalues of $G.$\\[1pt]
(i)\,\, If $\gl_{n-1}=\gl_{n}$ then $W$ determines the adjacency matrix of $G$ completely.\\
(ii)\, If $\gl_{n-1}\neq \gl_{n}$ then $W$ determines the adjacency matrix of at most two distinct graphs $G$ and $G^{*}$ with walk matrix $W.$
\end{thm}

{\sc Remarks 5.5:}  The exception in (ii) does occur, % for all $n$,
an example is given in  Appendix 8.3. There we obtained two different $8\times 8$ adjacency matrices with the same walk matrix of rank $6.$  In this particular example the corresponding graphs however are isomorphic. We conjecture that the two graphs in (ii) of the theorem are indeed always isomorphic. %This conjecture holds for the case of standard walk matrices, see Section~6.

{\sc Remarks 5.6:} The next cases to consider are graphs $G$ with walk matrix $W$ of rank $|V|-3$ or  $|V|-4.$  For ${\rm rank}(W)=|V|-3$ we have no example of a pair of adjacency matrices (whose graphs are isomorphic or otherwise) with the same walk matrix. For ${\rm rank}(W)=|V|-4$ we  have several pairs of non-isomorphic graphs of order  $9$ with the same walk matrix of rank $5.$ One example of this kind  is given in Appendix 8.4.

\medskip
{\it Proof of Theorem~\ref{n5.06}:}\,  By Theorem~\ref{n4.6}\, $W$ determines $M$ and hence the  main polynomial, say \begin{equation*}
  \main(x)=x^{n-2}+a_1x^{n-3}+\cdots+a_{n-3}x+a_{n-2}.
\end{equation*}
Let $\gl_{n-1}$ and $\gl_{n}$ be the two non-main eigenvalues and put
$$(x-\gl_{n-1})(x-\gl_n)=:x^2+b_1x+b_2.$$ Therefore, when $$\cha(x)=x^n+c_1x^{n-1}+c_2x^{n-2}+\cdots+c_{n-1}x+c_n\,,$$ then $\cha(x)=(x-\gl_{n-1})(x-\gl_n)\cdot \main(x)=(x^2+b_1x+b_2)\cdot \main(x)$ gives
\begin{equation*}
c_1=b_1+a_1=0 \text{\,\, and\,\,\,}
c_2=b_2+a_1b_1+a_2=-m,\\
\end{equation*}
where $m$ is the number of edges of $G.$ (We are using the well-known fact that $-c_{n-2}$ always is the number of edges of the graph, see~\cite{biggs, godsil2}.) Under the assumption (a) the column $W_{[1,1]}$ determines  the vertex degrees and from these we determine $m.$ Under assumption (b) this information is given anyhow.  Therefore $b_1=-a_1$ and $b_2=a_1^2-a_2-m$ and so  the two non-main eigenvalues \\[-10pt]
\begin{equation*}
\gl_{n-1,n}=\frac{a_1\pm\sqrt{4(a_2+m)-3a_1^2}}{2}.
\end{equation*}
of $G$ are determined by $W$ and $m.$

Next let
 $\f_{n-1}$ and $\f_{n}$ be non-main eigenvectors of $G$ for the eigenvalues $\gl_{n-1}$ and $\gl_{n}.$ Then $\f_{n-1}\cdot f_{n-1}^{\tt T}\,+\,\f_{n}\cdot f_{n}^{\tt T}$ is determined by $W$ according to \,(\ref{n5.14}). If $\gl_{n-1}=\gl_{n}$ we can determine $A$ by Theorem~\ref{n5.04}.

It remains to consider the case $\gl_{n-1}\neq\gl_{n}.$ Here we write \begin{eqnarray}\label{n5.17}
\gl_{n-1}\f_{n-1}\cdot f_{n-1}^{\tt T}\,+\,\gl_{n}\f_{n}\cdot f_{n}^{\tt T}&=&\gl_{n-1}(\f_{n-1}\cdot f_{n-1}^{\tt T}\,+\,\f_{n}\cdot f_{n}^{\tt T})\nonumber\\&+&(\gl_{n}-\gl_{n-1})\f_{n}\cdot f_{n}^{\tt T}\nonumber\\
&=&\gl_{n-1}({\rm I}_{n}- W_{[0,n-3]}\cdot W^{\dagger})\nonumber\\&+&(\gl_{n}-\gl_{n-1})\f_{n}\cdot f_{n}^{\tt T}\,,
\end{eqnarray}
using \,(\ref{n5.14}). Furthermore, from \,(\ref{n5.15})\, we obtain
\begin{eqnarray}\label{n5.18}
A&=&W_{[1,r]}\cdot W^{\dagger}\,\,+\,\,\sum_{j=n-1}^{n}\gl_{j}(\f_{j}\cdot \f_{j}^{\tt T})\nonumber\\&=&W_{[1,r]}\cdot W^{\dagger}+\gl_{n-1}({\rm I}_{n}- W_{[0,n-3]}\cdot W^{\dagger})\,+\,(\gl_{n}-\gl_{n-1})\f_{n}\cdot f_{n}^{\tt T}\,.
\end{eqnarray} In this equation $W_{[1,r]}\cdot W^{\dagger}+\gl_{n-1}({\rm I}_{n}- W_{[0,n-3]}\cdot W^{\dagger})$ is known and so we let $w_{1},\dots,w_{n}$ be the entries on the diagonal of that matrix. Similarly, denote the diagonal entries of $\f_{n}\cdot \f_{n}^{\tt T}$ by $\f_{n,1}^{2},\dots,\f_{n,n}^{2}$ where $\f_{n}=(\f_{n,1},\dots,\f_{n,n})^{\tt T}.$ Since the diagonal entries of $A$ are zero we have \begin{eqnarray}\label{n5.19} -(\gl_{n}-\gl_{n-1})^{-1}w_{i}&=&\f_{n,i}^{2}\,\,\,\text{\,\,for $i=1,\dots,n.$}\end{eqnarray}

Recall that $\f_{n}$ is a unit vector in the kernel of $W^{\tt T}.$ This space has dimension $2$ and so $\f_{n}$ depends on one parameter $t,$ i.e. an angle of rotation. The equation \,(\ref{n5.19})\, consists of a system of quadratic equations for $t$ and therefore it has at most two solutions. Hence there are at most two options for $\f_{n}.$ Finally return to \,(\ref{n5.18})\, to determine $A.$  \dne

\section{\sc Walk Equivalence}

It is natural to study graphs which have the same walk matrix  for  suitably chosen sets of  vertices.   We have the following characterisation in terms of the restrictions of the adjacency matrix.

\medskip
\begin{thm}\label{n6.01} Let $G$ and $G^{*}$ be graphs on the vertex set $V$ and let  $S,\,S^{*}\subseteq V.$ Denote by $W$ and $W^{*}$ the corresponding walk matrices, with restrictions  $A_{W}$ and $A^{*}_{W^{*}}$  respectively.   Then the following are equivalent  \\[1pt]
(i)\,\,\, $W=W^{*},$ \\
(ii)\,\, $G$ and $G^{*}$ have the same main eigenvalues and main eigenvectors for $S$  and \\
\phantom{(iii)\,} $S^{*}$ respectively, and \\
(iii)\, $S=S^{*}$ and $A_{W}=A^{*}_{W^{*}}.$\\[5pt]
In particular, if $W$ and $W^{*}$ are the standard walk matrices for $G$ and $G^{*},$ respectively, then
$W=W^{*}$ if and only if $A_{W}=A^{*}_{W^{*}}.$
\end{thm}

\pf If (i) holds then the first column of $W=W^{*}$ determines $S=S^{*}$ and from  \,(\ref{n5.15})\, we have $A_{W}=A^{*}_{W^{*}}.$ Hence (i) implies (iii). Next (ii) implies (i) by Theorem~\ref{n4.3}.

It remains to show that  (iii) implies (ii). Let $r$ and $r^{*}$ be the rank of $W$ and $W^{*}$ respectively.  By Proposition~\ref{n5.041} the condition $A_{W}=A^{*}_{W^{*}}$ implies $r^{*}=r$ or $r^{*}=r+1,$ without loss of generality. First suppose that $r^{*}=r$ when the proposition implies  $\mu_{i}^{*}=\mu_{i}$ and $\e_{i}^{*}=c_{i}\e_{i}$ with certain coefficients $c_{i}\in \ka,$ for $1\leq i\leq r.$  Since $S=S^{*}$ we have $x^{S}=\sum_{i} \e_{i}=x^{S^{*}}=\sum_{i} \e^{*}_{i}=\sum_{i} c_{i}\e_{i}$ and hence $\e_{i}^{*}=\e_{i}$ for $1\leq i\leq r.$ Thus (ii) holds in this case. Next suppose that $r^{*}=r+1$ when $\mu_{i}^{*}=\mu_{i}$ and $\e_{i}^{*}=c_{i}\e_{i}$ with certain coefficients $c_{i}\in \ka,$ for $1\leq i\leq r,$ and $\e_{r^*}=c_{r^{*}}f_{r+1},$ say. Since $S=S^{*}$ we have $x^{S}=\sum_{i=1}^{r} \e_{i}=x^{S^{*}}=\sum_{i=1}^{r} \e^{*}_{i}+\e_{r^{*}}=\sum_{i} c_{i}\e_{i}+\e_{r^{*}}.$ This implies $\e_{r^{*}}=0,$ a contradiction. Hence  (iii) implies (ii).  \dne

Clearly, it is essential to be able to reorder the vertices of a graph and change all  associated matrices accordingly. Let  $V=\{v_{1},\dots,v_{n}\}$ be the vertices of $G$
and let ${\rm Sym}_{n}$ be the symmetric group on $\{1,2,\dots,n\}.$ For  $g\in  {\rm Sym}_{n}$ we denote $g\!:\,i\to i^{g}$ for  $i=1,\dots n.$ From this we obtain permutations of the  elements and subsets of $V$ by setting $g\!:\,v_{i}\to v_{i}^{g}:=v_{i^{g}}$ for  $i=1,\dots n$ and $g\!:\, T\to T^{g}:=\{v^{g}\,|\,v\in T\}$ for any $T\subseteq S.$

The same permutation  gives rise to a new adjacency relation $\sim_{g}$ and a new graph $G^{*}=G^{g}$ by setting $v^{g}\sim_{g} u^{g}$ if and only if $v\sim u$ in $G.$ In this way  $g$ defines an isomorphism $g\!:\,G\to G^{*},$ denoted $G\simeq G^{*}.$

With $g \in {\rm Sym}_{n}$  we associate the $n\times n$ permutation matrix $P=P(g).$ It has the property that   for all $1\leq i\leq n$ we have $v_{j}=v_{i}^{g}$ if and only if $\vv_{j}=P\cdot \vv_{i}.$ (Here $\vv_{i}$ is the characteristic vector of $v_{i},$ etc..)  It follows that $P\cdot\x^{ T}=\x^{(T^{g})}$ for all $T\subseteq V$ when $\x^{T}$ denotes the characteristic vector of $T.$  The proof of the next lemma is left to the reader.

\medskip
\begin{lem}\label{n6.02} Let $G$ and $G^{*}$ be graphs on the vertex sets $V$ and $V^{*}$ with walk matrices $W$ and $W^{*}$ defined for certain sets $S\subseteq V$ and $S^{*}\subseteq V^{*},$ respectively. Suppose that $g\in {\rm Sym}(n)$ is a permutation  for which $G^{g}=G^{*}$ and $S^{g}=S^{*}.$ If $P:=P(g)$ then\\[1pt]
(i)\,\,\,\, $\e_{j}^{*}=P\e_{j}$ for all $1\leq j\leq r^{*}=r,$  $\e^{*}=P\e$ and $W^{*}=P\cdot W,$\\
(ii)\,\,  $A^{*}=PAP^{\tt T},$ and \\
(iii)  $A^{*}_{W^{*}}=PA_{W}P^{\tt T}.$
\end{lem}

\medskip
{\sc Definition:}\, Let $G$ and $G^{*}$ be graphs on the vertex sets $V$ and $V^{*},$ respectively,  and let   $S\subseteq V$ and $S^{*}\subseteq V^{*}.$   Denote the corresponding walk matrices by $W$ and $W^{*}.$ Then $(G,S)$ is {\it walk equivalent} to $(G^{*},S^{*}),$ denoted $(G,S)\sim (G^{*},S^{*}),$ if  there is a permutation matrix $P$ such that $W^{*}=P\cdot W.$ If $(G,S)\sim (G^{*},S^{*})$ and $S=V$ then $G$ is {\it walk equivalent} to $G^{*},$ denoted $G\sim G^{*}.$

%{\small \Fe 16. I worked very hard to make this definition as easy as possible. There is no problem when $S=V,$ when it becomes a relation on graphs. But in the general case, when $S\subset V$ things will become more complicated if we want an equivalence relation. The problem is transitivity: if we have $(G,S)\sim(G^{*},S^{*})$ and $(G^{*},S^{*})\sim(G^{**},S^{**})$  then all is fine, we can conclude that $(G,S)\sim(G^{**},S^{**}).$ I think it needs to be a relation on pairs of graphs with distinguished subset. \redne}

{\sc Remarks 6.1:}  In the definition, if we have $W^{*}=P\cdot W,$ then $S^{*}$ is determined by  $\e^{*}=P\cdot \e,$ the first column of $W^{*}.$  From this it follows easily that  $\sim$ indeed is an equivalence relation. To be precise, it is  a relation for graphs with a  distinguished vertex set. It turns into an equivalence relation on graphs per se only when $S=V.$

{\sc Remarks 6.2:}  If $G$ is isomorphic to $G^{*}$ via the permutation $g$ then it follows from Lemma~\ref{n6.02}\,  that $(G,S)\sim (G^{*},S^{*})$ for any $S\subseteq V$ when we set $S^{*}=S^{g}.$  Conversely however,  walk equivalence does not imply isomorphism. For instance, if $G$ and $G^{*}$ are regular graphs with the same valency then  $G\sim G^{*},$ and we may not conclude that $G$ is isomorphic to $G^{*}.$ Furthermore, if $G\sim G^{*}$ then  we can not conclude much in general about the walk equivalence of pairs $(G,S),\, (G^{*},S^{*})$ with  $S\subset V$ and $S^{*}\subset V^{*}.$ There are  simple examples for $G\sim G^{*}$ but $(G,S)\not\sim(G^{*},S^{*})$ for all $|S|,\,|S^{*}|<|V|.$  (For instance, let $G$ be the union of two $3$-cycles and $G^{*}$ a $6$-cycle.)

%{\sc Remarks 6.3:} In order to decide  whether $(G,S)$ is walk equivalent to $(G^{*},S^{*})$ it is useful to transform the corresponding walk matrices into a standard format.

{\sc Lexicographical Order:}\, If we want to decide  whether $(G,S)$ is walk equivalent to $(G^{*},S^{*})$ it is useful to transform the corresponding walk matrices into a standard format.

Let $W=W^{S}$ be the walk matrix of the graph $G$ for the vertex set $S.$ Then there is a permutation $g\in {\rm Sym}(n)$ of the rows of $W,$ with corresponding permutation matrix $P=P(g),$ so  that the rows of $PW$ are in lexicographical order.  The matrix
$$\lx(W):=PW$$  is the {\it lex-form} of $W$ with {\it reordering matrix} $P.$ Evidently $P$ is unique if and only if the rows of $W$ are pairwise distinct. It is also clear that  $(G,S)\sim (G^{*},S^{*})$ if and only if $\lx(W^{S})=\lx(W^{S^{*}}).$

In order to keep track of the reordering it is useful to append  the `label vector' ${\rm L}:=(v_{1},\dots,v_{n})^{\tt T}$ as last column to $W,$ obtaining $$W^{\ddag}=\big[W\,\big|\,{\rm L}\big].$$  The matrix  $QW^{\ddag}=\lx(W^{\ddag})$ then is  the {\it vertex lex-form} of $W.$

\medskip
{\sc Example 6.1:}\,  Let $G$ be the graphin Figure~3.1 on $V=\{v_{1},v_{2},v_{3},v_{4}\},$ let $S=\{v_{3}\}$ and $W:=W^{S}.$ Then  $$
W=\left(
\begin{array}{cccc}
 0 & 0 & 1 & 1 \\
 0 & 1 & 1 & 4 \\
 1 & 0 & 2 & 2 \\
 0 & 1 & 1 & 3 \\
\end{array}
\right)  \quad \text{and} \quad
W^{\ddag}=\left(
\begin{array}{ccccc}
 0 & 0 & 1 & 1 &v_{1}\\
 0 & 1 & 1 & 4 &v_{2}\\
 1 & 0 & 2 & 2 &v_{3} \\
 0 & 1 & 1 & 3 &v_{4}\\
\end{array}
\right)\,.
$$
Its vertex lex-form is
$$\lx\,W^{\ddag}=\left(
\begin{array}{ccccc}
1 & 0 & 2 & 2 &v_{3} \\
0 & 1 & 1 & 4 &v_{2}\\
0 & 1 & 1 & 3 &v_{4}\\
0 & 0 & 1 & 1 &v_{1}\\
\end{array}
\right) \,
$$
from which it follows that the reordering  permutation is $(v_{3},v_{1},v_{4})(v_{2}).$

\bigskip
\begin{thm}\label{n6.03} Let  $G$ and $G^{*}$ be graphs and suppose that there is a subset $S$ of vertices of $G$ such that $W^{S}$ has rank $\geq |V|-1.$ Then $G$ is isomorphic to $G^{*}$ if and only if $(G,S)\sim(G^{*},S^{*})$  for some vertex set $S^{*}$ of $G^{*}.$ Furthermore, if $G$ is isomorphic to $G^{*}$ then the isomorphism is determined uniquely from the lex forms of $W^{S}$ and $W^{S^{*}},$ unless $W^{S}$ has a repeated row when there are two isomorphisms.
\end{thm}

\pf If $G$ is isomorphic to $G^{*}$ via the map $g\in {\rm Sym}(n)$ the result follows from Lemma~\ref{n6.02} and Remark 6.2 above, taking $S^{*}=S^{g}.$ Conversely,
let $g\in {\rm Sym}(n)$ such that $P(g)W=W^{*}$ where $W=W^{S}$ and $W^{*}=W^{S^{*}}.$ Consider the graph $H=G^{g}.$ By Lemma~\ref{n6.02}  its  walk matrix is $P(g)W=W^{*}.$ Hence $H$ and $G^{*}$ have the same walk matrix. As the rank of $W^{*}$ is $\geq n-1$ it follows from Theorems~\ref{n5.01} and~\ref{n5.05}\, that $H$ and $G^{*}$ have the same adjacency matrix. Hence $H=G^{*}$ and so $G$ is isomorphic to $G^{*}.$ Suppose that $P(h)W=\lx(W)=P(h^{*})W^{*}$ for the corresponding reordering permutations $h$ and $h^{*}.$ It follows that $P(h)=P(h^{*})P(g)$ and so $g=(h^{*})^{-1}h.$ If the rows of $W$ (and hence $W^{*})$ are pairwise distinct then $h$ and $h^{*}$ are unique. %, see Remark~6.3 above.
If $W$ has repeated  rows then there is exactly one pair of repeated rows since  ${\rm rank}(W)\geq n-1$ by assumption. In this case  $g=(h^{*})^{-1}h$ is unique up to the transposition interchanging these rows. \dne.

The following is a special case of Theorem~\ref{n6.03}:

\bigskip
\begin{thm}\label{n6.04} Let  $G$ be a graph on the vertex set $V$ and let $S,\,S^{*}$ be subsets of $V.$ Suppose that $W^{S}$ has rank $\geq |V|-1.$ Then there is an automorphism $g$ of $G$ with  $S^{g}=S^{*}$ if and only if $\lx(W^{S})=\lx(W^{S^{*}}).$ \end{thm}

%\Fe 9. This is quite a nice result. But you had additional ideas about graph isomorphisms? Could we add any conjecture, maybe prove it? \dner

{\sc Example 6.2:}\,  The graph in Figure~3.1 has the following walk matrices for $S=\{3\}$ and  $S=\{4\},$ respectively,
{\small $$
W^{\{3\}}=\left(
\begin{array}{cccc}
 0 & 0 & 1 & 1 \\
 0 & 1 & 1 & 4 \\
 1 & 0 & 2 & 2 \\
 0 & 1 & 1 & 3 \\
\end{array}
\right)
\,\quad \text{and}\,\quad
W^{\{4\}}=\left(
\begin{array}{cccc}
 0 & 0 & 1 & 1 \\
 0 & 1 & 1 & 4 \\
 0 & 1 & 1 & 3 \\
 1 & 0 & 2 & 2 \\
\end{array}
\right).$$}
\hspace{-2mm} These matrices are lex-equivalent and have rank $\geq 3.$ Therefore Theorem~\ref{n6.04} implies that there is an automorphism of $G$  which interchanges $3$ and $4.$

{\it Proof of Theorem~\ref{n1.2}:}  This  follows from Theorem~\ref{n6.03}. \dne

%{\it Proofs of Theorem~\ref{n1.2} and \ref{n1.3}}: CHECK WHEN WE HAVE REWRITTEN SECTION 1. These results follow from Theorems~\ref{n6.03} and \ref{n5.02}. \dne

\section{Random Graphs}

Let $\mathcal{P}$ be a property of finite undirected simple graphs $G$. Then we say that $\mathcal{P}$ holds {\it almost always}, or that {\it almost all graphs have property $\mathcal{P}$} if the probability for $\mathcal{P}$ to hold, as a function of the order $n$ of $G,$ tends to $1$ as $n$ tends to infinity. The following result is due to O'Rourke and Touri\cite{orourke1} and is based on Tao and Vu\cite{tao}.  Recall that the standard walk  matrix of $G$ with vertex set $V$ is $W^{V}:$

\medskip
\begin{thm}\label{DD1} For almost all graphs the standard walk matrix is invertible.
\end{thm}

From Theorems~\ref{n5.01} and ~\ref{n6.04} we have immediately the  following two conclusions. Recall that $G$ and $G^{*}$ are walk-equivalent if their standard walk matrices have the same lex form.

\medskip
\begin{thm}\label{DD2} (i) For almost all graphs $G$  the standard walk matrix of $G$ determines the  adjacency matrix of $G.$ \\[5pt]
(ii) For almost all graphs $G$ we have that  $G$ is isomorphic to the  graph $G^{*}$ if and only if $G$ is walk-equivalent to $G^{*}.$ \end{thm}

{\it Proof of Theorem~\ref{n1.3}}: The result follows from Theorem~\ref{DD2}. \dne

Let $G$ be a graph with characteristic polynomial ${\rm}char_{G}(x).$ In Theorem~\ref{n5.01a} we have noted that if ${\rm}char_{G}(x)$ is irreducible then $W^{S}$ is invertible for any vertex set $S,$ and so Theorem~\ref{DD1} holds more generally.
In the literature there are several paper in which this irreducibility problem is considered from a probabilistic point of view, see~\cite{Vu2, diet, orourke2 , vuwood}. In fact, there is the

\medskip
\begin{con}\label{CJa} For  most all graphs the characteristic polynomial is irreducible.
 \end{con}

The characteristic polynomial of a graph is irreducible if and only if all  eigenvalues are simple and the Galois group of the graph acts transitively on these eigenvalues.  For recent papers on the Galois group of integers polynomials see Dietmann~\cite{diet}. We have the following stronger

\medskip
\begin{con}\label{CJb} For almost all graphs $G$ we have $\Gal(G)={\rm Sym}(n)$
 where $n$ is the order of $G.$
 \end{con}

 This conjecture is supported by a theorem of Van der Waerden~\cite{vdw} which states  that for almost all monic integer polynomials the Galois group is isomorphic to ${\rm Sym}(n),$ where $n$ is the degree of the polynomial. There are examples of graphs with irreducible characteristic polynomial where $\Gal(G)\neq {\rm Sym}(n).$ But these are difficult to find, in most practical  computations the group turns out to be  ${\rm Sym}(n).$

 \bigskip
\begin{thm}\label{DD3} If the Conjecture~\ref{CJa} holds then the following is true for almost all graphs $G$ and arbitrary vertex set $S$ of $G:$ For any graph $G^{*}$ there is an isomorphism $G\to G^{*}$ if and only if $(G,S)\sim(G^{*},S^{*})$  for some vertex set $S^{*}$ of $G^{*}.$ Furthermore, if $G$ is isomorphic to $G^{*}$ then the isomorphism can determined from the lex forms of $W^{S}$ and $W^{S^{*}}.$
\end{thm}

%%%%%%%%%%%%%%%%%%%%%%%%%%%%%%%%%%%%%%%%%%%%%%%%%%%%%%
%%%%%%%%%%%%%   APPENDIX ITEMS  BEGIN %%%%%%%%%%%%%%%%%%%%%%%%%%%

\section{\sc Appendix}

We give details and examples of graphs and walk matrices with particular features mentioned in earlier sections of the paper.

%\begin{ex}\label{ex-same-E&Eigenspace}
{\sc Appendix 8.1:} \,The following are two non-isomorphic graphs $G$ and $G^{*}$  on $8$ vertices with the same eigenspaces.  Let $G$ and $G^{*}$ be the graphs shown in Fig. 8.1 with adjacency matrices $A$ and $A^{*}.$ Let
{\small \begin{equation*}
  T=\begin{bmatrix}
 1 & -1 & 0 & 0 & 0 & -1 & 0 & -1 \\
 1 & -1 & 0 & 0 & -1 & 0 & 0 & 1 \\
 1 & -1 & 0 & 0 & 0 & 1 & 0 & -1 \\
 1 & -1 & 0 & 0 & 1 & 0 & 0 & 1 \\
 1 & 1 & 0 & -1 & 0 & 0 & -1 & 0 \\
 1 & 1 & -1 & 0 & 0 & 0 & 1 & 0 \\
 1 & 1 & 0 & 1 & 0 & 0 & -1 & 0 \\
 1 & 1 & 1 & 0 & 0 & 0 & 1 & 0 \\
\end{bmatrix}.
\end{equation*}}\\[8pt]
It is straightforward to verify
that $T^{-1}AT=\text{diag}(6,-2,0,0,0,0,-2,-2)$ and $T^{-1}A^{*}T=\text{diag}(5,-3,-1,-1,-1,-1,1,1).$ Therefore the columns of $T$ are simultaneous eigenvectors of $G$ and $G^{*}. $ Hence $G$ and $G^{*}$ are eigenspace equivalent but not isomorphic.
\begin{figure}[!h]
  \centering
  %\includegraphics[width=12cm]{SameEiv}\\
    \includegraphics[width=12cm]{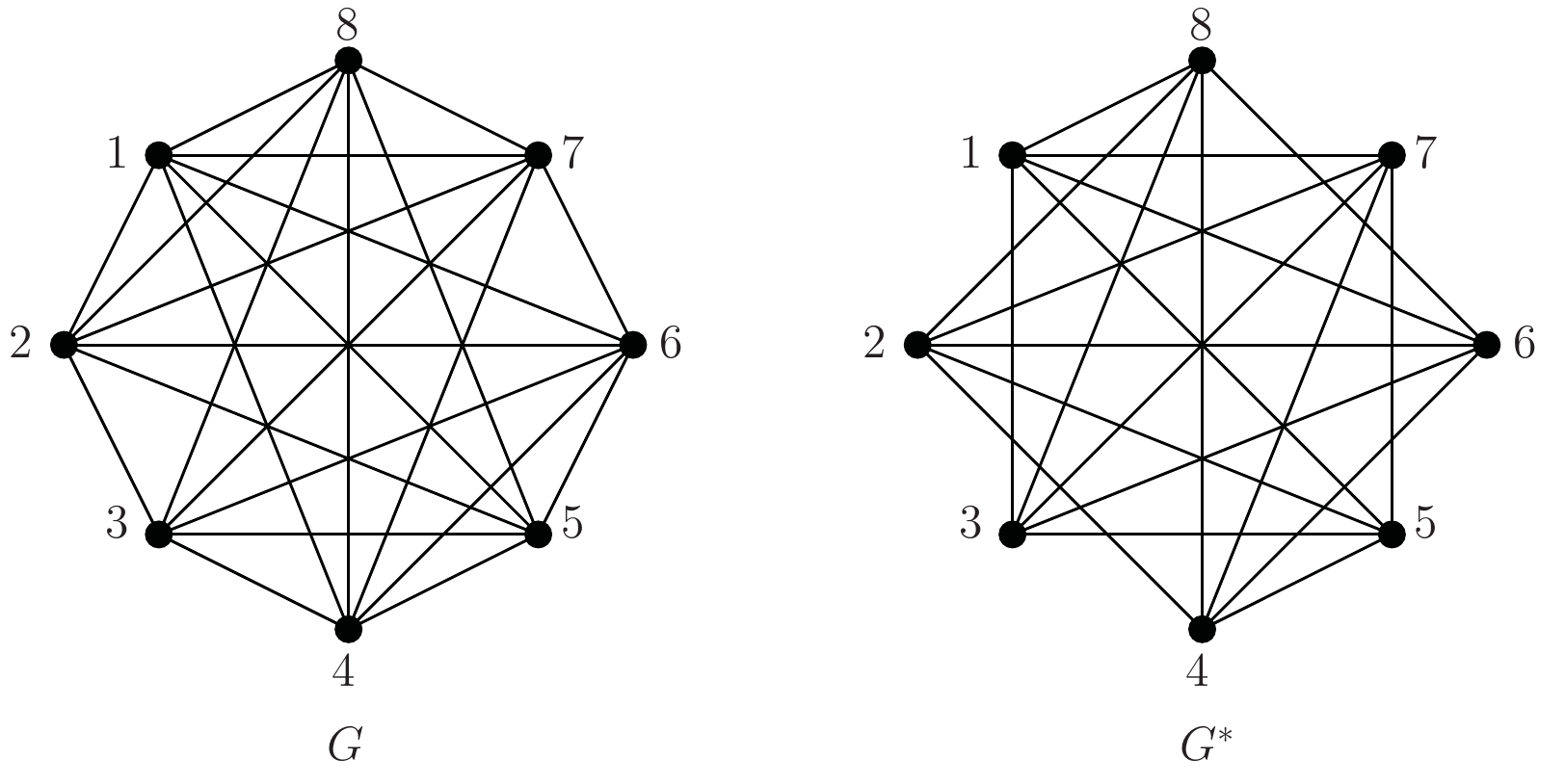}\\
  \sc Fig 8.1: Graphs with the Same Eigenspaces
  %\caption{\sc Graphs with the Same Eigenspaces}%\label{fig-SE1}
\end{figure}\\[20pt]

%\begin{figure}[!h]
 % \centering
 % \includegraphics[width=5cm]{SameEiv1}\qquad \qquad\includegraphics[width=5cm]{SameEiv2}\\
 % \caption{\sc Graphs with the Same Eigenspaces}%\label{fig-SE1}
%\end{figure}

\bigskip \bigskip \bigskip \bigskip

{\sc Appendix 8.2:} \, The following is an example of a pair of non-isomorphic graphs on $8$ vertices which have the same main eigenvectors for $S=V$ but different main eigenvalues. It  appears in L. Collins and I. Sciriha~\cite{Sciriha}.

Let $G$ and $G^{*}$ be as shown in Fig. 8.2, and let
$$e_1=(\frac{1+\sqrt{5}}{2}, \frac{1+\sqrt{5}}{2}, \frac{1+\sqrt{5}}{2},\frac{1+\sqrt{5}}{2},1,1,1,1)^{\tt T},$$$$e_2=( \frac{1-\sqrt{5}}{2}, \frac{1-\sqrt{5}}{2}, \frac{1-\sqrt{5}}{2},\frac{1-\sqrt{5}}{2},1,1,1,1)^{\tt T}.$$
It is straightforward to verify that
$e_1$ and $e_2$ are the main eigenvectors of $G$ and $G^{*}$ with two main eigenvalues
$$1+\sqrt{5},1-\sqrt{5}\,\,\,\,\text{and}\,\,\,\,\frac{3}{2}(1+\sqrt{5}),\frac{3}{2}(1-\sqrt{5}),$$ respectively.

%\newpage
\begin{figure}[!h]
  \centering
  % Requires \usepackage{graphicx}    %  \includegraphics[width=12cm]{SameMainEiv}
  \includegraphics[width=12cm]{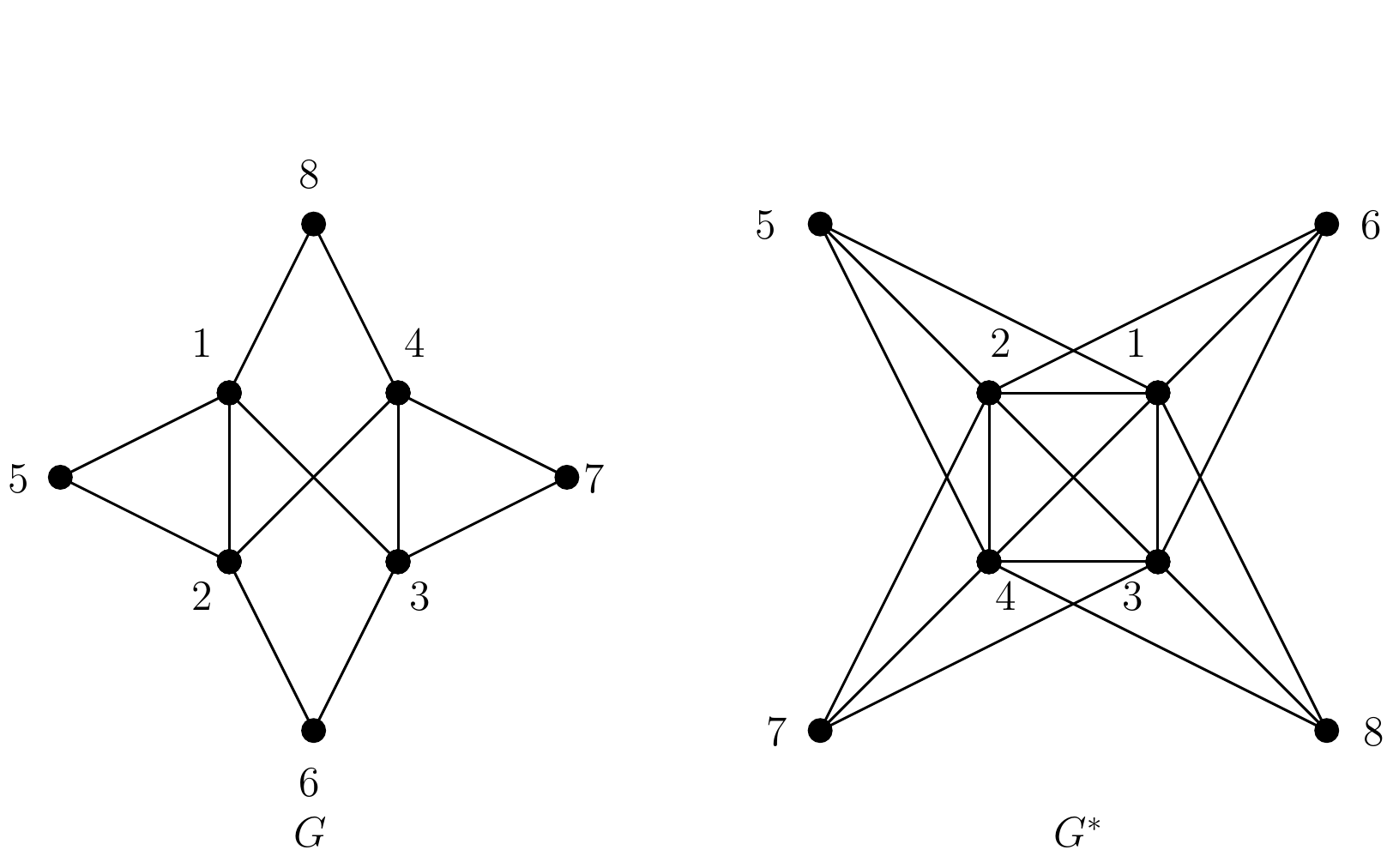}\\ {\sc Fig 8.2: Graphs with same Main Eigenvectors but\\ different Main Eigenvalues}
  %\caption{Graphs with same main eigenvectors but different main eigenvalues}\label{}
\end{figure}

%\begin{figure}[!h]
%  \centering
 % \includegraphics[width=6cm]{SameMain}\qquad \includegraphics[width=5cm]{SameEiv4}\\
 % \caption{Graphs with same main eigenvectors but different main eigenvalues}\label{}
%\end{figure}

\bigskip\bigskip
{\sc Appendix 8.3:} The following two  adjacency matrices on $8$ vertices \\
\begin{equation*}
  A_1=\begin{bmatrix}
 0 & 0 & 0 & 1 & 1 & 0 & 1 & 0 \\
 0 & 0 & 1 & 0 & 1 & 0 & 0 & 1 \\
 0 & 1 & 0 & 0 & 1 & 1 & 0 & 0 \\
 1 & 0 & 0 & 0 & 0 & 1 & 0 & 0 \\
 1 & 1 & 1 & 0 & 0 & 0 & 1 & 0 \\
 0 & 0 & 1 & 1 & 0 & 0 & 0 & 0 \\
 1 & 0 & 0 & 0 & 1 & 0 & 0 & 0 \\
 0 & 1 & 0 & 0 & 0 & 0 & 0 & 0 \\
\end{bmatrix},\
 A_2=\begin{bmatrix}
 0 & 0 & 1 & 0 & 1 & 0 & 0 & 1 \\
 0 & 0 & 0 & 1 & 1 & 0 & 1 & 0 \\
 1 & 0 & 0 & 0 & 1 & 1 & 0 & 0 \\
 0 & 1 & 0 & 0 & 0 & 1 & 0 & 0 \\
 1 & 1 & 1 & 0 & 0 & 0 & 1 & 0 \\
 0 & 0 & 1 & 1 & 0 & 0 & 0 & 0 \\
 0 & 1 & 0 & 0 & 1 & 0 & 0 & 0 \\
 1 & 0 & 0 & 0 & 0 & 0 & 0 & 0 \\
 \end{bmatrix},
\end{equation*}\\
give rise to the same  standard walk matrix $W$ of rank $6,$\\
\begin{equation*}
  W=\begin{bmatrix}
 1 & 3 & 8 & 23 & 64 & 181 & 506 & 1425 \\
 1 & 3 & 8 & 23 & 64 & 181 & 506 & 1425 \\
 1 & 3 & 9 & 24 & 69 & 190 & 539 & 1502 \\
 1 & 2 & 5 & 13 & 37 & 101 & 287 & 797 \\
 1 & 4 & 11 & 32 & 89 & 252 & 705 & 1984 \\
 1 & 2 & 5 & 14 & 37 & 106 & 291 & 826 \\
 1 & 2 & 7 & 19 & 55 & 153 & 433 & 1211 \\
 1 & 1 & 3 & 8 & 23 & 64 & 181 & 506 \\
 \end{bmatrix}.
\end{equation*}

In this case the graphs for $A$ and $A^{*}$ are isomorphic to each other. This  graph is drawn in Figure 8.3.
\begin{figure}[!h]
  \centering
  \includegraphics[width=6cm]{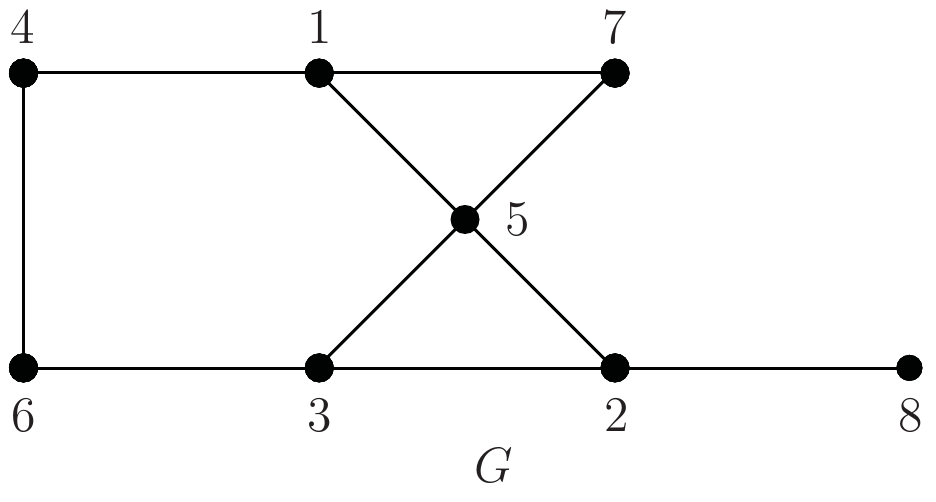}\\{\sc Fig 8.3: Two Different Adjacency Matrices with \\the  same standard Walk Matrix}
  %\caption{Two Different Adjacency Matrices with the  Same Walk Matrix}\label{fig-v8SameW1}
\end{figure}\\

\medskip\bigskip\bigskip
{\sc Appendix 8.4:}\, The smallest non-isomorphic graphs with the same standard walk matrix are the graphs $G$ and $G^*$ in Figure 8.4a on $7$ vertices. Their standard walk matrix is
$$
 W=W^*=\begin{bmatrix}
 1 & 4 & 11 & 35 & 104 & 318 & 960 \\
 1 & 3 & 9 & 27 & 82 & 248 & 752 \\
 1 & 2 & 7 & 20 & 62 & 186 & 566 \\
 1 & 2 & 8 & 22 & 70 & 208 & 636 \\
 1 & 2 & 7 & 20 & 62 & 186 & 566 \\
 1 & 3 & 9 & 27 & 82 & 248 & 752 \\
 1 & 4 & 11 & 35 & 104 & 318 & 960 \\
\end{bmatrix}.$$

\begin{figure}[!h]
  \centering
  \includegraphics[width=12cm]{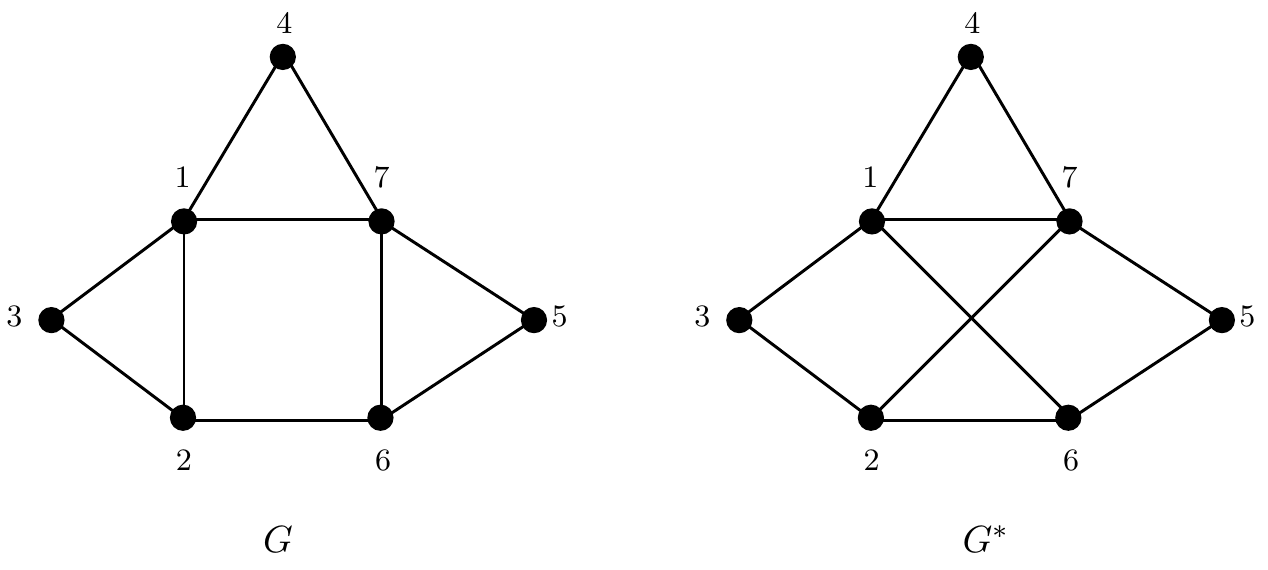}\\[8pt]

  {\sc Fig. 8.4{\rm a}: A pair of non-isomorphic Graphs with the same standard \\Walk Matrix of ${\rm rank}(W)=n-3=4.$}
  %\caption{Graphs with the same Walk Matrix and ${\rm rank}(W)=n-3$}\label{fig-SameW}
\end{figure}

\smallskip
Another example are the following non-isomorphic graphs $G$ and $G^{*}=\overline{G}$ on $n=9$ vertices in Figure 8.4b. Their  standard walk matrix\\
\begin{equation*}
  W=W^{*}=\begin{bmatrix}
 1 & 4 & 18 & 72 & 300 & 1222 & 5028 & 20586 & 84480 \\
 1 & 4 & 16 & 67 & 272 & 1121 & 4586 & 18827 & 77162 \\
 1 & 5 & 20 & 83 & 339 & 1393 & 5707 & 23413 & 95989 \\
 1 & 5 & 20 & 83 & 339 & 1393 & 5707 & 23413 & 95989 \\
 1 & 4 & 16 & 67 & 272 & 1121 & 4586 & 18827 & 77162 \\
 1 & 3 & 13 & 52 & 215 & 878 & 3607 & 14778 & 60625 \\
 1 & 4 & 16 & 65 & 267 & 1093 & 4485 & 18385 & 75403 \\
 1 & 4 & 16 & 65 & 267 & 1093 & 4485 & 18385 & 75403 \\
 1 & 3 & 13 & 52 & 215 & 878 & 3607 & 14778 & 60625 \\
  \end{bmatrix}.
\end{equation*}\\
has rank $n-4=5,$ see also  \cite{ Godsil-Mckay}.\\

\begin{figure}[!h]
  \centering
  %  \includegraphics[width=12cm]{NonIsoSameW}
  \includegraphics[width=12cm]{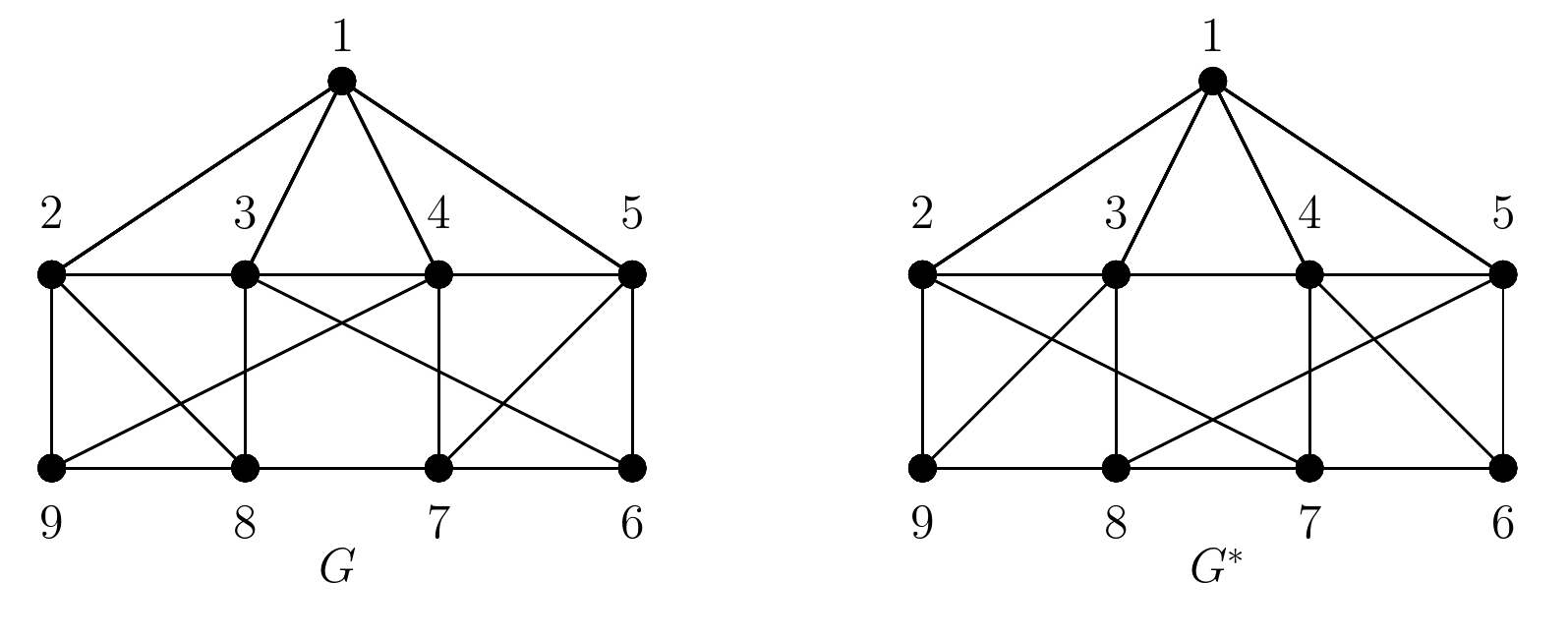} \\{\sc Fig. 8.4{\rm b}: A Pair of Non-isomorphic Graphs \\with the same Standard Walk Matrix of ${\rm rank}(W)=n-4=5.$}
  %\caption{A pair of non-isomorphic graphs with the same walk matrix}\label{fig-v9nonisosamew}
\end{figure}

%\bigskip\bigskip\bigskip\bigskip

\bigskip\bigskip
{\sc Appendix 8.5:} We return to the graph $G$ in Figure~3.1 in order to demonstrate in an example how the Galois group and the automorphism group act on the main eigenvectors and main eigenvalues of the graph.  The  characteristic polynomial of $G$ is $(1 + x) (x^3-x^{2}-3x+1),$ it has the roots $$\gl_{0}=-1,\, \,\gl_{1}=-1.48...,\,\, \gl_{2}=0.31...,\,\,\text{and}\,\,\,\gl_{3}=2.17...$$
with  corresponding eigenvectors \\[5pt]
{\small $$a_{0}=\left(\begin{array}{c}0 \\0 \\-1 \\1\end{array}\right),\quad
a_{1}=\left(\begin{array}{c}-0.67... \\1 \\-0.40... \\-0.40...\end{array}\right),\quad
a_{2}=\left(\begin{array}{c}3.21... \\1 \\-1.45... \\-1.45...\end{array}\right)\,\,\text{and}\,\,\,\,
a_{3}=\left(\begin{array}{c}0.46... \\1 \\0.85... \\0.85...\end{array}\right).$$}\\[5pt]
Observe that the graph automorphism  which fixes $v_{1},\,v_{2}$ and interchanges $v_{3},\,v_{4}$ maps $a_{0}$ to $-a_{0}$ and fixes $a_{1},\,a_{2},\,a_{3}.$ (Clearly,  graph automorphisms always fix all  eigenspaces: In matrix terms, if $P$ is a permutation matrix with $A\,P=P\,A$ and if $Ax=\lambda x$ then $A\, Px=P\, Ax=\lambda Px.$) The Galois group \,${\rm Gal}(G)$\, of $(x^3-x^{2}-3x+1)$ is the symmetric group on $\{\gl_{1},\gl_{2},\,\gl_{3}\},$ its elements permute ${\rm Eig}(A,\gl_{1}),\,{\rm Eig}(A,\gl_{2}),\,{\rm Eig}(A,\gl_{3})$ in the same way.

Next we find the main eigenvectors for $S=\{v_{1}\},$ for instance. The main polynomial of this set is   $(x^3-x^{2}-3x+1).$ We write $\e=(1,0,0,0)^{\tt T}$ in terms of its main eigenvectors, $$ \e=(1,0,0,0)^{\tt T}=\e_{1}+\e_{2}+\e_{3}=c_{1}a_{1}+c_{2}a_{2}+c_{3}a_{3}$$ for certain coefficients $c_{i}\in \ka,$ the splitting field of $ (x^3-x^{2}-3x+1).$ In fact, $c_{1}=-0.37..,$ $c_{2}=0.20...$ and $c_{3}=0.17...$\,. By Lemma~\ref{NN311}\, we have that  ${\rm Gal}(G)$ permutes the $\e_{i}$ and so we obtain the explicit action of ${\rm Gal}(G)$ on the  main eigenvectors of  $S=\{v_{1}\}.$

\bigskip
{\sc Appendix 8.6:} By computation we found that the characteristic polynomial of pairs of cospectral graphs of order $n<8$ is always reducible. The following in an example for $n=8$ with a pair of non-isomorphic graphs  $G,$ $G^{*}$ with  irreducible characteristic polynomial
$$x^8-10x^6-4 x^5+24 x^4+8x^3-16x^2+1.$$
Their complements have the same characteristic polynomial
$$x^8-18 x^6-26 x^5+26 x^4+42 x^3-16 x^2-16 x+4$$ and this polynomial is  again irreducible.

\begin{figure}[!h]
  \centering
  \includegraphics[width=10cm]{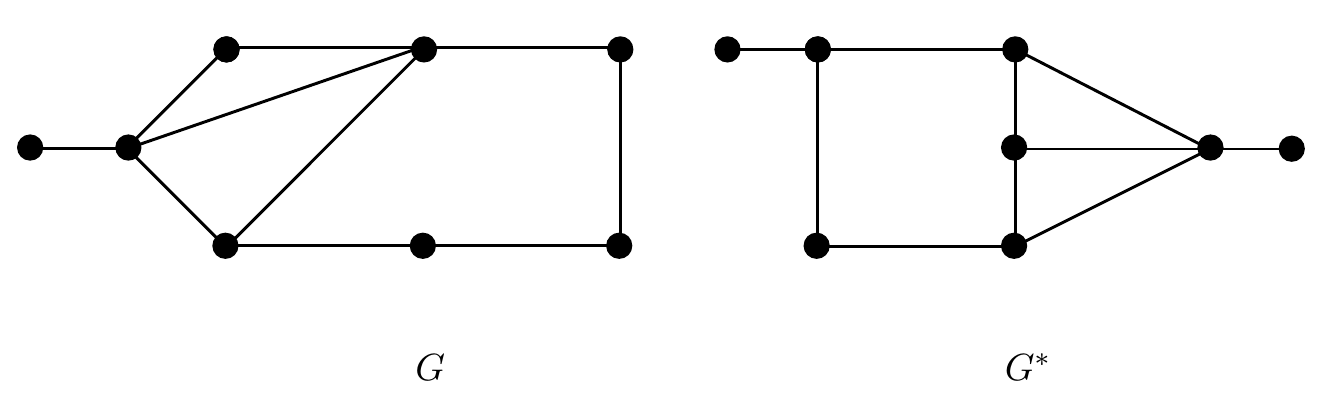}\\{\sc Fig. 8.6: Two non-isomorphic Graphs of order $6$ \\with the same Irreducible Characteristic Polynomial}

%\caption{Cospectrals Graphs with the same Irreducible Characteristic Polynomial }\label{fig-SameIrPoly}
\end{figure}

% \Fe Last Question, should we add an acknowledgment like the following, mentioning your students who help with the computations? \redne

 {\sc Acknowledgment:} We thank Professor Wei Wang for helpful comments. Several examples of graphs in the Appendix where obtained by Zenan Du. The first author is thankful for the support from the National Natural Science Foundation of China, the Fundamental Research Funds for the Central Universities and the Foundation of China Scholarship Council.

%%%%%%%%%%%%%%%%%%%%%%%%%%%%%%%%%%%%%%%%%%%%%%%%%%%%%%
%%%%%%%%%%%%%   APPENDIX ITEMS  END %%%%%%%%%%%%%%%%%%%%%%%%%%%

\end{document}